\documentclass[4pt, oneside]{article}
\usepackage{amsmath,amssymb,amsthm,units,stmaryrd,upgreek}

\usepackage{yfonts}
\usepackage{units,stmaryrd}
\usepackage{color}
\usepackage{graphicx}
\usepackage[all]{xy}

\newtheorem{definition}{Definition}[section]
\newtheorem{lemm}{Lemma}[section]
\newtheorem{prop}{Proposition}[section]
\newtheorem{cor}{Corollary}[section]

\newtheorem{theorem}{Theorem}[section]
\newtheorem{examp}{Example}[section]

\def\lp{\llparenthesis}
\def\rp{\rrparenthesis}
\def\lb{\left\llbracket}
\def\rb{\right\rrbracket}
\def\reductive{{r}}

\newcommand{\glp}{{\ensuremath{\mathsf{GLP}}}\xspace}

\usepackage{xspace}

\def\le{{\ell}}
\def\ex{e}

\def\fmodels{\xymatrix{
\ar@{|=}[r]^{<\omega}&
}
}
\def\nmodels{\xymatrix{
\ar@{|=}[r]^{N}&
}
}
\def\<{\left <}

\def\nc{{\Box}}
\def\ps{{\Diamond}}

\def\>{\right >}

\def\cbra{\left \{}
\def\cket{\right \}}

\DeclareSymbolFont{AMSb}{U}{msb}{m}{n}
\DeclareMathSymbol{\N}{\mathbin}{AMSb}{"4E}
\DeclareMathSymbol{\Z}{\mathbin}{AMSb}{"5A}
\DeclareMathSymbol{\R}{\mathbin}{AMSb}{"52}
\DeclareMathSymbol{\Q}{\mathbin}{AMSb}{"51}
\DeclareMathSymbol{\I}{\mathbin}{AMSb}{"49}
\DeclareMathSymbol{\C}{\mathbin}{AMSb}{"43}

\newcommand{\comment}[2]{{\color{red}{\bf #2}}}

\begin{document}

\title{The polytopologies of transfinite provability logic
      }

\author{David Fern\'{a}ndez-Duque\footnote{Group for Computational Logic, Universidad de Sevilla, dfduque@us.es}}

\maketitle

\begin{abstract}
Provability logics are modal or polymodal systems designed for modeling the behavior of G\"odel's provability predicate and its natural extensions. If $\Lambda$ is any ordinal, the G\"odel-L\"ob calculus $\mathsf{GLP}_\Lambda$ contains one modality $[\lambda]$ for each $\lambda<\Lambda$, representing provability predicates of increasing strength. $\mathsf{GLP}_\omega$ has no Kripke models, but it is sound and complete for its topological semantics, as was shown by Icard for the variable-free fragment and more recently by Beklemishev and Gabelaia for the full logic.

In this paper we generalize Beklemishev and Gabelaia's result to ${\sf GLP}_\Lambda$ for countable $\Lambda$. We also introduce {\em provability ambiances}, which are topological models where valuations of formulas are restricted. With this we show completeness of ${\sf GLP}_\Lambda$ for the class of provability ambiances based on Icard polytopologies.
\end{abstract}


\section{Introduction}

Provability logic interprets modal operators as provability predicates in order to study the structure of formal theories, reading the modal formula $\nc\phi$ as {\em the theory $T$ proves $\phi$}. In \cite{Solovay:1976}, Solovay proved that if $T$ is able to do a reasonable amount of arithmetic, the set of validities over the unimodal language is given by the G\"odel-L\"ob logic $\sf GL$, written $\glp_1$ in the current paper's notation. This logic may also be interpreted over {\em scattered spaces} (where every non-empty subset has an isolated point), thus giving provability a surprising connection to topology. However, in practice these semantics are somewhat heavy-handed for such a logic, which already has finite Kripke models based on transitive, well-founded frames \cite{Segerberg:1971}. 

For Japaridze's polymodal provability logic, the story is not as simple. It is an extension of $\sf GL$ known as $\glp$ or, in our notation, $\glp_\omega$ \cite{Japaridze:1988}. Here one considers countably many provability modalities $[n]$, for $n<\omega$. The formula {$[n]\phi$} could be interpreted (for example) as {\em $\phi$ is derivable using $\omega$-rules of depth at most $n$}. There is great interest in $\mathsf{GLP}$ since these logics are quite powerful and useful; Beklemishev has shown how $\glp$ can be used to perform ordinal analysis of Peano Arithmetic and its natural subtheories \cite{Beklemishev:2004}.

However, the logic is no longer as easy to work with as in the unimodal case. As we shall discuss later, it has no non-trivial Kripke frames. Thus the topological interpretation of the logic gives a reasonable alternative, but even then we do not get an immediate solution to the problem. In fact, the existence of so-called {\em canonical ordinal models} for these theories goes well beyond $\sf ZFC$, as shown by Blass \cite{blass}, Beklemishev \cite{glb} and in recent unpublished work by Bagaria.

There are, however, polytopologies based on ordinals for which $\glp=\glp_\omega$ is sound and complete, as shown by Beklemishev and Gabelaia \cite{BeklGabel:2011}. The proof of this difficult result requires some heavy machinery including Zorn's lemma, so the resulting spaces are non-constructive. There are also simpler spaces which provide semantics for the {\em closed fragment}, where no free variables occur; these were introduced by Icard \cite{Icard:2009} and are closely tied to Ignatiev's Kripke model for the same fragment \cite{Ignatiev:1993}.

Our goal is to show how the constructions from \cite{BeklGabel:2011} may be extended to the logics $\glp_\Lambda$, where $\Lambda$ is an arbitrary ordinal.  Here, one has transfinitely many provability operators, which as in the case of $\glp_\omega$ represent derivability in stronger and stronger theories. Indeed, Beklemishev and Gabelaia's techniques carry over smoothly to the transfinite setting, and rather than give a new, self-contained completeness proof, we shall state the necessary results from \cite{BeklGabel:2011} without proof in order to focus on applying these techniques beyond $\omega$. A key point is the computation of the higher-order rank functions, which give us upper and lower bounds on the ordinals we need in order to build models. We shall also show how the use of non-constructive topologies may be circumvented and replaced by Icard topologies by passing to a more general class of models called {\em ambiances}.

\paragraph{Layout.} In Section \ref{syntax} we give a quick overview of the logics $\mathsf{GLP}_\Lambda$, and Section \ref{topsem} reviews topological semantics. Section \ref{ordinalintro} then states some basic facts about ordinal arithmetic that we shall need.

Section \ref{ranks} introduces the most important functions in the study of GLP-spaces, {\em ranks} and {\em $d$-maps}. Then, Section \ref{secic} discusses Icard ambiances and Section \ref{secsimp} simple ambiances, the minimal structures in our framework.

After this, Section \ref{bgsec} discusses{ Beklemishev-Gabelaia spaces}, which are particularly well-behaved GLP-spaces. In Section \ref{secred}, we discuss and construct reductive functions, an important type of $d$-map, and Section \ref{secop} establishes a series of operations on ambiances which are used for constructing models.

We then go on to review the logic $\sf J$ in Section \ref{jsec}, which is a key ingredient in the completeness proof presented in Section \ref{seccomp}. Finally, Section \ref{lowsec} uses {\em worms,} which are special variable-free formulas related to ordinals, to give a lower bound on the rank of models.


\section{The logic $\mathsf{GLP}_\Lambda$}\label{syntax}

Given any ordinal $\Lambda$, we can define a provability logic with modalities in $\Lambda$. Formulas of the language $\mathsf{L}_\Lambda$ are built from $\top$ and a countable set of propositional variables $\mathbb P$ using Boolean connectives $\neg,\wedge,\vee,\to$ and a modality $[\xi]$ for each $\xi<\Lambda$. As is customary, we use $\<\xi\>$ as a shorthand for $\neg[\xi]\neg$.

The logic $\mathsf{GLP}_\Lambda$ is then given by the following rules and axioms:
\begin{enumerate}
\item all propositional tautologies{,}
\item $[\xi](\phi\to\psi)\to([\xi]\phi\to[\xi]\psi)$ for all $\xi<\Lambda${,}
\item {$[\xi]([\xi]\phi \to \phi)\to[\xi]\phi$ for all $\xi<\Lambda$}{,}
\item $[\xi]\phi\to[\zeta]\phi$ for $\xi<\zeta<\Lambda${,}\label{glpfour}
\item $\<\xi\>\phi\to [\zeta]\<\xi\>\phi$ for $\xi<\zeta<\Lambda$,
\item modus ponens and
\item necessitation for each $[\xi]$.
\end{enumerate}

Note that the unimodal ${\sf GLP}_1$ is the standard G\"odel-L\"ob logic $\sf GL$. Let us write ${\rm sub}(\phi)$ for the set of subformulas of $\phi$. Then, we say that $\lambda$ {\em appears} in $\phi$ if there is some formula $\psi$ such that $[\lambda]\psi\in{\rm sub}(\phi)$. It is evident that only finitely many ordinals may appear in any formula $\phi$; sometimes it is convenient to ignore all other ordinals. To this end we define the {\em condensation} of $\phi$ as follows:

\begin{definition}
Given a formula $\phi\in{\sf L}_\Lambda$ such that
\[\lambda_0<\lambda_1<\hdots<\lambda_{N-1}\]
are the ordinals appearing in $\phi$, we define a formula $\phi^{\sf c}$ (the {\em condensation} of $\phi$) as the result of replacing every operator $[\lambda_n]$ in $\phi$ by $[n]$.
\end{definition}

As it turns out, the formula $\phi^{\sf c}$ is derivable if and only if $\phi$ is. One direction, which we will not need in this paper, is non-trivial and proven in \cite{paper0}; the other is quite straightforward and will be used later.

\begin{lemm}\label{condenselemm}
If $\phi$ is a formula such that there are $N$ ordinals appearing in $\phi$ then $\glp_N\vdash\phi^{\sf c}$ implies that $\glp_\Lambda\vdash\phi$.
\end{lemm}

This fact may be proven by uniformly substituting $[\lambda_n]$ for $[n]$ in a derivation of $\phi^{\sf c}$; we omit the details. Condensations will allow us to focus only on `relevant' ordinals when analyzing formulas.

We shall also work with Kripke semantics. A {\em Kripke frame} is a structure $\mathfrak F=\<W,\<R_n\>_{n<N}\>$, where $W$ is a set and $\<R_n\>_{n<N}$ a family of binary relations on $W$. A {\em valuation} on $\mathfrak F$ is a function $\lb\cdot\rb:{\sf L}_\Lambda\to \mathcal P(W)$ such that
\[
\begin{array}{lcl}
\lb\bot\rb&=&\varnothing\\\\
\lb\neg\phi\rb&=&W\setminus\lb\phi\rb\\\\
\lb\phi\wedge\psi\rb&=&\comment{W}{}\lb\phi\rb\cap\lb\psi\rb\\\\
\lb\<n\>\phi\rb&=&R^{-1}_n\lb\phi\rb.
\end{array}
\]

A {\em Kripke model} is a Kripke frame equipped with a valuation $\lb\cdot\rb$. Note that propositional variables may be assigned arbitrary subsets of $W$. If $\mathfrak M=\langle\mathfrak F,\lb\cdot\rb\rangle$ is a model, we may write $\<\mathfrak M,x\>\models\psi$ instead of $x\in\lb \psi\rb$. As usual, $\phi$ is {\em satisfied} on $\mathfrak M$ if $\lb\phi\rb\not=\varnothing$, and {\em true} on $\mathfrak M$ if $\lb\phi\rb=W$. It is {\em valid} on a frame $\mathfrak F$ if it is true on every model based on $\mathfrak F$.

It is well-known that L\"ob's axiom is valid on $\mathfrak F$ whenever $R^{-1}_n$ is well-founded and transitive \cite{Segerberg:1971}, in which case we denote it by $<_n$. However, constructing models of $\glp_\Lambda$ is substantially more difficult than constructing models of $\mathsf{GL}$; the full logic $\mathsf{GLP}_\Lambda$ cannot be sound and complete with respect to any class of Kripke frames. Indeed, let $\mathfrak F=\langle W,\<<_\xi\>_{\xi<\lambda}\rangle$ be a polymodal frame.

Then, it is not too hard to check that
\begin{enumerate}
\item \label{FrameCondOne}L\"ob's axiom $[\xi]([\xi]\phi \to \phi)\to [\xi]\phi$ is valid if and only if $<_\xi$ is well-founded and transitive{,}

\item \label{FrameCondTwo}the axiom $[\xi] \phi \to[\zeta] \phi$ {for $\xi\leq\zeta$} is valid if and only if, whenever $w<_{\zeta}v${}, then $w<_{\xi}v$, and

\item \label{FrameCondThree}$\langle \xi \rangle \phi \to [\zeta] \langle \xi \rangle \phi$ {for $\xi<\zeta$} is valid if, whenever $v<_{\zeta}w $, $u<_{\xi}w$ and $\xi<\zeta$, then $u<_{\xi}v$.
\end{enumerate}

Suppose that for $\xi<\zeta$, there are two worlds such that $w<_\zeta v$. Then from \ref{FrameCondTwo} we see that $w<_\xi v$, while from \ref{FrameCondThree} this implies that $w<_\xi w$. But this clearly violates \ref{FrameCondOne}. Hence if $\mathfrak F\models \mathsf{GLP}$, it follows that all accessibility relations (except possibly $<_0$) are empty.

This observation makes the topological completeness of $\glp_\omega$ established in \cite{BeklGabel:2011} particularly surprising. Moreover, as we shall see, the techniques introduced there readily extend to the transfinite. To show this, let us begin by reviewing the topological semantics of provability logic.


\section{Topological semantics}\label{topsem}

Recall that a {\em topological space} is a pair $\mathfrak X=\<X,\mathcal T\>$ where $\mathcal T\subseteq \mathcal P(X)$ is a family of sets called `open' such that
\begin{enumerate}
\item $\varnothing,X\in\mathcal T$
\item if $U,V\in\mathcal T,$ then $U\cap V\in \mathcal T$ and
\item if $\mathcal U\subseteq\mathcal T$ then $\bigcup\mathcal U\in\mathcal T$.
\end{enumerate}

Given $A\subseteq X$ and $x\in A$, we say $x$ is a {\em limit point} of $A$ if, for all $U\in\mathcal T$ such that $x\in U$, we have that $(A\setminus \cbra x\cket)\cap U\not=\varnothing$. We denote the set of limit points of $A$ by ${d} A$, and call it the `derived set' of $A$. We can define topological semantics for modal logic by interpreting Boolean operators in the usual way and setting
\[\lb\ps\psi\rb_\mathfrak X={d}\lb\psi\rb_\mathfrak X.\]

In order to interpret provability logic, we will need to consider {scattered} spaces. A topological space $\langle X,\mathcal T\rangle$ is {\em scattered} if every non-empty subset $A$ of $X$ has an isolated point; that is, there exist $x\in A$ and a neighborhood $U$ of $x$ (i.e., $x\in U\in\mathcal T$) such that $U\cap A=\cbra x\cket$.

Many interesting examples of scattered spaces come from ordinals. The simplest is the {\em initial segment topology}. If $\Theta$ is an ordinal, we use $\Theta_{0}$ to denote the structure $\langle\Theta,\mathcal T\rangle$, where $\mathcal T$ consists of all downward-closed subsets of $\Theta$. It is very easy to check that $\Theta_{0}$ is a scattered topological space, for if $A\subseteq \Theta$ is non-empty, then the least element of $A$ is isolated in $A$.

A second important example is the {\em interval topology}. This is generated by all intervals on $\Theta$ of the form $[0,\beta]$ or $(\alpha,\beta]$. The interval topology extends the initial segment topology, and it is straightforward to check that if $\mathcal T$ is scattered and $\mathcal T'$ is any refinement of $\mathcal T$ (i.e., $\mathcal T\subseteq\mathcal T'$), then $\mathcal T'$ is scattered as well. We will denote $\Theta$ equipped with the interval topology by $\Theta_1$.

Now, in order to interpret $\glp_\Lambda$ for $\Lambda>1$, we need to consider {\em polytopological spaces}. A {polytopological space} is a structure $\mathfrak X=\<X,\<\mathcal T_\lambda\>_{\lambda<\Lambda}\>,$ where $\Lambda$ is an ordinal and each $\mathcal T_\lambda$ is a topology. The derived set operator corresponding to $\mathcal T_\lambda$ shall be denoted $d_\lambda$. We may also write $\mathfrak X_\lambda$ instead of $\langle X,\mathcal T_\lambda\rangle$.

There are Kripke-incomplete modal logics which nevertheless are complete for {\em general Kripke frames}, which are Kripke frames where valuations are restricted to a special algebra of sets. A similar idea will prove useful in order to give constructive semantics of ${\sf GLP}_\Lambda$. The following definition describes the algebras we shall use:

\begin{definition}[$d$-algebra]\
A {\em $d$-algebra} over a polytopological space $\mathfrak X=\langle X,\langle\mathcal T_\lambda\rangle_{\lambda<\Lambda}\rangle$ is a collection of sets $\mathcal A\subseteq \mathcal P(X)$ which form a Boolean algebra under the standard set-theoretic operations and such that, whenever $\lambda<\Lambda$ and $S\in\mathcal A$ it follows that $d_\lambda S\in\mathcal A$.
\end{definition}

Below we introduce ambiances, which will be the basis of our semantics; they are a slight generalizarion of polytopological models, which correspond to the special case where $\mathcal A=\mathcal P(X)$.

\begin{definition}[Ambiance]
An {\em ambiance} is a structure
\[\mathfrak X=\langle X,\vec{\mathcal T},\mathcal A\rangle\]
consisting of a polytopological space equipped with a $d$-algebra $\mathcal A$.

\end{definition}

If $\vec{\mathcal T}=\langle\mathcal T_\lambda\rangle_{\lambda<\Lambda}$, we may also say $\mathfrak X$ is a {\em $\Lambda$-ambiance}. The operator $d_\lambda$ will be used to interpret $\langle \lambda\rangle$:

\begin{definition}
Let  $\mathfrak X=\langle X,\vec{\cal T},\mathcal A\rangle$ be a $\Lambda$-ambiance.

A {\em valuation} on $\mathfrak X$ is a function $\lb\cdot\rb:{\sf L}_\Lambda\to \mathcal A$ defined as in the case of Kripke semantics except that
\[\lb \langle\lambda\rangle\phi\rb=d_\lambda\lb\phi\rb.\]

A {\em polytopological model} is an ambiance equipped with a valuation.
\end{definition}

Let us check the conditions under which $\glp_\Lambda$ is sound for a given polytopological space.
\begin{lemm}\label{glpspace}
Let $\Lambda$ be an ordinal and $\mathfrak X=\langle X,\vec{\cal T},\mathcal A\rangle$ be an ambiance.

Then,
\begin{enumerate}
\item \label{CorrespondenceOne} L\"ob's axiom $[\xi]([\xi]\phi \to \phi)\to [\xi]\phi$ is valid on $\mathfrak X$ whenever $\<X,\mathcal T_\xi\>$ is scattered{,}

\item \label{CorrespondenceTwo} the axiom $[\xi] \phi \to [\zeta] \phi$ {for $\xi\leq\zeta$} is valid whenever $\mathcal T_\xi\subseteq \mathcal T_\zeta$ and

\item \label{CorrespondenceThree} $\langle \xi \rangle \phi \to [\zeta] \langle \xi \rangle \phi$ {for $\xi<\zeta$} is valid if $d_\xi A\in \mathcal T_\zeta$ whenever $A\in\mathcal A$.
\end{enumerate}

\end{lemm}
\proof
See, for example, \cite{BeklGabel:2011}.
\endproof

An ambiance satisfying the above properties will be called a {\em provability ambiance}.

When referring to topologies, we use the words {\em extension} and {\em refinement} indistinctly. We may also speak of {refinements} of {\em spaces} rather than refinements of topologies: $\langle X',\mathcal T'\rangle$ is a refinement of $\langle X,\mathcal T\rangle$ if $X=X'$ and $\mathcal T\subseteq\mathcal T'$. Thus the condition for $[\xi] \phi \to [\zeta] \phi$ can be rewritten as {``$\mathcal T_\zeta$ is a refinement of $\mathcal T_\xi$''}.

Conditions \ref{CorrespondenceTwo} and \ref{CorrespondenceThree} suggest a very natural candidate for $\mathcal T_{\xi+1}$ whenever $\mathcal T_\xi$ is given; namely, the {\em least} topology that will satisfy all axioms.

\begin{definition}[$d^{\mathcal A}{\mathcal T}$]
Given a $\mathfrak X=\langle X,\mathcal T\rangle$ and a $d$-algebra $\mathcal A$ on $\mathfrak X$, we define $d^{\mathcal A}\mathcal T$ to be the topology on $X$ generated by
\[\mathcal T\cup\{dS:S\in\mathcal A\}.\]

We will denote $\langle X,d^{\mathcal A}\mathcal T\rangle$ by $d^{\mathcal A}\mathfrak X$.
\end{definition}

As in \cite{BeklGabel:2011}, we shall write $\mathfrak X^+$ instead of $d^{\mathcal P(X)}\mathfrak X$. The above definition suggests natural candidate topologies for $\mathcal T_\lambda$, at least for successor $\lambda$. For limit $\lambda$ we need to consider {\em joins} of topologies.

If $\vec{\mathcal T}=\mathcal T_{\xi<\lambda}$ is an increasing sequence of topologies, then $\mathcal U=\bigcup_{\xi<\lambda}\mathcal T$ is typically not a topology. Although it is always closed under finite unions and intersections, it need not be closed under arbitrary unions. However, $\vec{\mathcal T}$ does generate a least topology $\mathcal J=\bigsqcup_{\xi<\lambda}\mathcal T$ (its `join') containing all $\mathcal T_\xi$, by closing $\mathcal U$ under arbitrary unions. The elements of $\mathcal J$ are then of the form $\bigcup_{\lambda<\Lambda}U_\lambda$ with $U_\lambda\in\mathcal T_\lambda$. In other words, $\mathcal U$ forms a basis for $\mathcal J$, so that $O$ is open in $\mathcal J$ if and only if for every $x\in O$ there are $\lambda<\Lambda$ and $V\in\mathcal T_\lambda$ such that $x\in V\subseteq O$.

Due to the monotonicity axiom, this is the least topology we can choose at limit stages:

\begin{definition}\label{tminus}
Given an ambiance $\mathfrak X=\langle X,\vec{\mathcal T},\mathcal A\rangle$ and an ordinal $\xi>0$, define $\mathcal T^-_\xi$ to be
\begin{itemize}
\item $d^{\mathcal A}\mathcal T_\zeta$ if $\xi=\zeta+1$\\
\item $\displaystyle\bigsqcup_{\zeta<\xi}\mathcal T_\zeta$ if $\xi$ is a limit ordinal.
\end{itemize}
\end{definition}
A na\"ive strategy for building models of $\sf GLP$ consists of always choosing the least possible topology at each stage; a structure ${\mathfrak X}=\langle\Theta,\langle\mathcal T_\lambda\rangle_{\lambda<\Lambda}\rangle$ is a {\em canonical ordinal model} if ${\mathfrak X}_0=\Theta_{1}$ and for all $\lambda<\Lambda$, $\mathcal T_{\lambda}=\mathcal T^-_\lambda$. In a canonical ordinal model the topology $\mathcal T_{\xi+1}=\mathcal T_\xi^+$ is usually much bigger than $\mathcal T_\xi$, since we are adding many new closed sets as opens. For example:
\begin{lemm}
Given an ordinal $\Theta,$ $\Theta_{0}^+=\Theta_1$.
\end{lemm}

The reader may wish to prove this directly as an exercise; we will not give such a proof as it is a special case of Lemma \ref{dotchar}. After this the topologies increase very quickly; if $\Theta$ is any countable ordinal and $\mathcal T$ is the interval topology, then $\mathcal T^+$ is discrete. Moreover, the question of whether $\glp_2$ is complete for its class of canonical ordinal models is independent of $\rm ZFC$ \cite{blass,glb}. In recent unpublished work, Bagaria has characterized non-trivial ordinal models for ${\sf GLP}_n$ in terms of large cardinals.

Thus, making the topologies as small as possible at each step is not the best strategy, so it is convenient to consider other alternatives. In \cite{BeklGabel:2011}, Beklemishev and Gabelaia realized the highly unintuitive fact that if we make each topology {\em as large as possible} then subsequent topologies become much smaller! Thus they obtain spaces where $\mathcal T_n\supsetneq \mathcal T^-_n$ for each $n>1$. As we shall see, this idea readily extends beyond $\glp_\omega$; perhaps the most technically challenging aspect of such an extension lies in the new computations with ordinals that arise.


\section{Operations on ordinals}\label{ordinalintro}

Before continuing, let us give a brief review of some notions of ordinal arithmetic as well as some useful functions in the study of provability logic. We skip most proofs; for further details on ordinal arithmetic, we refer the reader to a text such as \cite{PohlersBook}, while the material on hyperexponentials and hyperlogarithms is treated in detail in \cite{hyperexp}.

We assume familiarity with ordinal sums, products and exponents. We shall also use the following operations:

\begin{lemm}\label{theorem:BasicPropertiesOrdinalArithmetic}\

\begin{enumerate}
\item Whenever $\zeta {<} \xi$, there exists a unique ordinal $\eta$ such that $\zeta + \eta = \xi$. We will denote this unique $\eta$ by $-\zeta + \xi$.

\item
Given $\xi>0$, there exist ordinals $\alpha,\beta$ such that $\xi = \alpha + \omega^{\beta}$. The value of $\beta$ is uniquely defined.
We will denote this unique $\beta$ by $\le \xi$.
\end{enumerate}
\end{lemm}

In previous work my colleague Joost Joosten and I realized that there were some particularly useful functions that arise when studying provability logics. They are {\em hyperexponentials} and {\em hyperlogarithms}, and are a form of transfinite iteration of the functions $-1+\omega^\xi$ {and $\le$, respectively}. These iterations have been used in \cite{WellOrdersI} for describing well-orders in the Japardize algebra and in \cite{glpmodels} for defining models of the variable-free fragment of $\glp_\Lambda$. They will be essential in defining our semantics. We give only a very brief overview, but \cite{hyperexp} gives a thorough and detailed presentation.

We shall denote the class of all ordinals by $\mathsf{On}$ and the class of limit ordinals by $\mathsf{Lim}$.

\begin{definition}\label{hyperdef}
Let $\ex(\xi)=-1+\omega^\xi$. Then, we define the {\em hyperexponentials} $\langle\ex^\zeta\rangle_{\zeta\in\sf On}$ as the unique family of normal\footnote{That is, strictly increasing and continuous.} functions such that
\begin{enumerate}
\item $\ex^1=\ex$
\item $\ex^{\alpha+\beta}=\ex^\alpha \ex^\beta$ for all ordinals $\alpha,\beta$
\item $\langle\ex^\zeta\rangle_{\zeta\in\sf On}$ is pointwise minimal amongst all families of normal functions satisfying the above clauses\footnote{That is, if $\langle g^\zeta\rangle_{\zeta\in\sf On}$ is a family of functions satisfying conditions 1 and 2, then for all ordinals $\xi,\zeta$, $\ex^\zeta\xi\leq g^\zeta\xi$.}.
\end{enumerate}
\end{definition}

It is not obvious that such a family of functions exists, but a detailed construction is given in \cite{hyperexp}, where the following is also proven:

\begin{prop}[Properties of {hyperexponentials}]\label{prophyp}
The family of functions $\langle\ex^\xi\rangle_{\xi\in\mathsf{On}}$ has the following properties:
\begin{enumerate}
\item $e^0$ is the identity,\label{prophypone}
\item $e^\xi 0=0$ for all $\xi$,\label{prophyptwo}
\item given $\xi\in\mathsf{On}$ and $\lambda\in \sf Lim$, $\ex^{\xi}\lambda=\displaystyle \lim_{\eta\to\lambda}\ex^\xi \eta$ and\label{prophypthree}
\item if $\lambda\in\sf Lim$ and $\vartheta\in\sf On$, $e^\lambda(\vartheta+1)=\lim_{\eta\to\lambda}e^\eta(e^\lambda(\vartheta)+1)$.\label{prophypfive}
\end{enumerate}
\end{prop}

\begin{examp}\label{ExEpsilon}
We have that $e(0)=-1+\omega^0=-1+1=0$ and $e(1)=-1+\omega^1=\omega$. Then, $e^21=ee1=e(\omega)=-1+\omega^\omega=\omega^\omega$, and continuing in this fashion one sees that
\[e^n1=\underbrace{\omega^{\omega^{\cdot^{\cdot^{\cdot^\omega}}}}}_n.\]
We know from Proposition \ref{prophyp}.\ref{prophyptwo} that $e^\omega 0=0$, and in view of Proposition \ref{prophyp}.\ref{prophypfive},
\[e^\omega 1=\lim_{n\to\omega}e^n(e^\omega(0)+1)= \lim_{n\to\omega}e^n1 =\lim_{n\to\omega} \underbrace{\omega^{\omega^{\cdot^{\cdot^{\cdot^\omega}}}}}_n,\]
usually denoted $\varepsilon_0$. Meanwhile
\[e^\omega 2=\lim_{n\to\omega}e^n(\varepsilon_0+1)=\lim_{n\to\omega} \underbrace{\omega^{\cdot^{\cdot^{\cdot^{\omega^{\varepsilon_0+1}}}}}}_n=\varepsilon_1,\] and more generally $\varepsilon_n=e^\omega(n+1)$.

Finally, by Proposition \ref{prophyp}.\ref{prophypthree} we see that \[e^\omega\omega=\lim_{n\to\omega}e^\omega n=\lim_{n\to \omega}\varepsilon_n=\varepsilon_\omega.\]
This can be generalized to obtain $\varepsilon_\xi=e^\omega(1+\xi)$ for every ordinal $\xi$.
\end{examp}

Closely related to hyperexponentials are hyperlogarithms. Below, an {\em initial function} is one mapping initial segments to initial segments.

\begin{definition}[Hyperlogarithms]
We define the sequence $\langle \le^{\xi}\rangle_{\xi\in\mathsf{On}}$ to be the unique family of initial functions such that
\begin{enumerate}
\item $\le^1=\le$,
\item $\le^{\alpha+\beta}=\le^{\beta}\le^{\alpha}$ for all ordinals $\alpha,\beta$,
\item $\langle \le^{\xi}\rangle_{\xi\in\mathsf{On}}$ is pointwise maximal among all families of functions satisfying the above clauses.
\end{enumerate}
\end{definition}

The following properties of hyperlogarithms will be used throughout the text and are not too difficult to check:

\begin{prop}\label{LemmLogProp}
The hyperlogarithms $\langle\le^\xi\rangle_{\xi\in {\sf On}}$ have the following properties:
\begin{enumerate}
\item $\le^0$ is the identity,
\item  If $\alpha,\delta>0$ and $\gamma$ is any ordinal, then $\le^\alpha(\gamma+\delta)=\le^\alpha\delta$.\label{LemmLogPropOne}
\item For any ordinal $\gamma$, the sequence $\langle\le^\xi\gamma \rangle_{\xi\in{\sf On}}$ is non-increasing.
\end{enumerate}\label{LemmLogPropTwo}
\end{prop}

Observe that hyperexponentials are typically not surjective, hence not right-invertible. However, they are injective, thus left-invertible, and hyperlogarithms provide particularly well-behaved left inverses.

\begin{lemm}\label{leftinv}
If $\xi<\zeta$, then $\le^\xi\ex^\zeta=\ex^{-\xi+\zeta}$ and $\le^\zeta\ex^\xi=\le^{-\xi+\zeta}$.

Further, whenever $\alpha<e^\xi\beta$, it follows that $\le^\xi\alpha<\beta$.
\end{lemm}

We may also use the contrapositive form of the above, that is, whenever $\beta\leq \le^\xi\alpha$, then $e^\xi\alpha\leq \beta$. Note that it also follows from this that when $\beta < \le^\xi\alpha$, then $e^\xi\alpha< \beta$, for if we had $e^\xi\alpha= \beta$ then also $\le^\xi e^\xi\alpha=\alpha= \le^\xi\beta$.

\begin{examp}
Let us compute the sequence $\langle \le^\xi\gamma\rangle_{\xi\in \sf On}$ for $\gamma=\varepsilon_{\omega\cdot 3}+\varepsilon_{\omega\cdot 2}$. First we have that $\le^0\gamma=\gamma$ since $\le^0$ is the identity.

The $\varepsilon$-numbers have the property that they are fixed under the map $\xi\mapsto\omega^\xi$, so we may also write $\gamma$ as $\varepsilon_{\omega\cdot 3}+\omega^{\varepsilon_{\omega\cdot 2}}$. In view of this, $\le^1 \gamma=\le\gamma= \varepsilon_{\omega\cdot 2}$. After this, $\le^2\gamma=\le\le \gamma=\le \varepsilon_{\omega\cdot 2}=\varepsilon_{\omega\cdot 2}$, and continuing inductively we see that $\le^n\gamma=\varepsilon_{\omega\cdot 2}$ for all $n<\omega$.

To go beyond $\omega$ and in view of Example \ref{ExEpsilon}, we may write $\gamma$ as $e^\omega(\omega\cdot 3)+e^\omega(\omega\cdot 2)$, so that $\le \gamma=e^\omega(\omega\cdot 2)$. Therefore, by Lemma \ref{leftinv}, $\le^\omega\gamma=\le^\omega e^\omega(\omega\cdot 2)=\omega\cdot 2$. Then, $\le^{\omega+1}\gamma=\le\le^\omega\gamma= 1$, since $\omega\cdot 2=\omega+\omega^1$, and thus $\le^{\omega+2}\gamma=\le\le^{\omega+1}\gamma=\le 1=0$, since $1=\omega^0$. From here on we obtain $\le^\xi\gamma= 0$ for all $\xi>\omega+1$.

In summary, the sequence $\langle\le^\xi\gamma\rangle_{\xi\in \sf On}$ has the following form:
\[\underbrace{\varepsilon_{\omega\cdot 3}+\varepsilon_{\omega\cdot 2}}_{\xi=0},\underbrace{\varepsilon_{\omega\cdot 2},\varepsilon_{\omega\cdot 2},\hdots}_{1\leq \xi<\omega},
\underbrace{\omega\cdot 2}_{\xi=\omega},
\underbrace{1}_{\xi=\omega+1},
\underbrace{0,0,\hdots}_{\xi>\omega+1}
\]
\end{examp}

There is a close relation between the iterates $\ex^{\omega^\gamma}\xi$ and Veblen functions; this is also described in detail in \cite{hyperexp}. For example:

\begin{lemm}\label{rangelemm}
An ordinal $\xi$ lies in the range of $\ex^{\omega^\gamma}$ if and only if, for all $\delta<\gamma$, we have that $\xi=\ex^{\omega^\delta}\xi$. In particular, $\ex^{\omega^{\gamma+1}}$ enumerates the fixpoints of $\ex^{\omega^{\gamma}}$.
\end{lemm}

This has some consequences which will prove to be very useful to us:

\begin{lemm}\label{LemmLogConverge}
If $\Lambda=\alpha+\omega^\beta$ is a limit ordinal and $\xi$ is any ordinal, then there exists $\lambda<\Lambda$ such that $\le^\vartheta\xi=e^{\omega^\beta}\le^\Lambda\xi$ for all $\vartheta\in[\lambda,\Lambda).$
\end{lemm}

\proof
Since $\le^\vartheta\xi$ is non-increasing on $\vartheta$, there must be some $\lambda<\Lambda$ such that $\le^\vartheta\xi=\le^\lambda\xi$ for all $\vartheta\in[\lambda,\Lambda)$; clearly we may pick $\lambda>\alpha$. Observe then that for all $\delta<\beta$ we have that \[\le^{\omega^\delta}\le^\lambda\xi= \le^{\lambda+\omega^\delta}\xi =\le^\lambda\xi,\]
so that by Lemma \ref{leftinv}, $\le^\lambda\xi\geq e^{\omega^\delta}\le^\lambda\xi$; since $e^{\omega^\delta}$ is normal this means that $\le^\lambda\xi = e^{\omega^\delta}\le^\lambda\xi$. By Lemma \ref{rangelemm} we have that $\ell^\lambda\xi=e^{\omega^\beta}\eta$ for some $\eta>0$, and applying $\le^{\omega^\beta}$ on both sides we obtain
\[\le^\Lambda\xi= \le^{\lambda+\omega^\beta}\xi =\le^{\omega^\beta}\le^\lambda\xi =
\eta.\]
Thus for $\vartheta\in[\lambda,\Lambda)$ we have that \[\le^\vartheta\xi=\le^\lambda\xi=e^{\omega^\beta}\eta=
e^{\omega^\beta}\le^\Lambda\xi,\]
as needed.
\endproof

\begin{lemm}\label{falls}
Suppose that for ordinals $\vartheta,\gamma$ and additively indecomposable $\Lambda$ we have that $\vartheta\in(e^\Lambda\gamma,e^\Lambda(\gamma+1))$.

Then, there exists $\lambda<\Lambda$ such that $\le^\lambda\vartheta\leq e^\Lambda\gamma$.
\end{lemm}

\proof
By Lemma \ref{LemmLogConverge} we have that $\le^\lambda\vartheta=e^\Lambda\le^\Lambda\vartheta$ for some $\lambda<\Lambda$. Since $\le^\lambda\vartheta\leq\vartheta< e^\Lambda(\gamma+1)$ and $e^\Lambda$ is normal, we must have $\le^\Lambda\vartheta\leq \gamma$ and thus $\le^\lambda\vartheta\leq e^\Lambda\gamma$.
\endproof

To conclude this section, let us discuss {\em simple functions.} We will often be faced with families of inequalities of the form $\{\le^{\alpha_n}\xi>\beta_n\}_{n<N}$, and need to describe the ordinals $\xi$ satisfying such constraints. Simple functions will be used to gather such inequalities into a single object. They will play a crucial role throughout the paper, as they provide a convenient, flexible tool for solving many of the problems that will arise later.

\begin{definition}
A {\em simple function} is a partial function $s\colon\Lambda\dashrightarrow\Theta$ with finite domain, where $\Theta,\Lambda$ are ordinals. We denote the domain of $s$ by ${\rm dom}(s)$.
\end{definition}

If $r,s$ are simple functions, we define $r\sqcup s$ to be the simple function with domain ${\rm dom}(r)\cup {\rm dom}(s)$ given by
\[r\sqcup s(\lambda)=
\begin{cases}
r(\lambda)&\text{if $\lambda\in{\rm dom}(r)\setminus {\rm dom}(s)$,}\\
s(\lambda)&\text{if $\lambda\in{\rm dom}(s)\setminus {\rm dom}(r)$,}\\
\max\{r(\lambda),s(\lambda)\}&\text{if $\lambda\in{\rm dom}(r)\cap {\rm dom}(s)$.}
\end{cases}\]

We write $s\sqsubseteq\alpha$ if for $\lambda=\max({\rm dom}(s))$ we have that $s(\lambda)\leq \le^{\lambda}{\alpha}$ and for all $\xi\in{\rm dom}(s)\setminus\{\lambda\}$, $s(\xi)<\le^\xi \alpha$. If moreover $s(\lambda)< \le^{\lambda}{\alpha}$, we instead write $s\sqsubset\alpha$. If $\mathop{\rm dom}(s)=\varnothing$, we also set $s\sqsubseteq \alpha$ and $s\sqsubset\alpha$ for all $\alpha$. Given a simple function $s$, we define $\lceil s\rceil$ to be the least ordinal $\sigma$ such that $s\sqsubseteq \sigma$. Note that if $\mathop{\rm dom}(s)=\varnothing,$ then $\lceil s\rceil=0.$

The following lemma originally appeared in \cite{glpmodels} in a different presentation.

\begin{lemm}\label{ceilbound}
Let $r$ be a simple function with non-empty domain and $\lambda=\max({\rm dom}(r))$.

Then,
\begin{enumerate}
\item
the ordinal $\lceil r\rceil$ is defined and $\le^\lambda\lceil r\rceil=r(\lambda)$, and\label{ceilbound1}
\item whenever $r\sqsubseteq\alpha$ and $\xi\leq\lambda$, we have that $\le^\xi\lceil r\rceil\leq \le^\xi\alpha$.\label{ceilbound2}
\end{enumerate}
\end{lemm}

\proof
It will be convenient for our proof to define $\mathop{\rm dom}^+(r)=\{0\}\cup \mathop{\rm dom}(r)$. We will proceed to construct an ordinal $\vartheta$ such that $r\sqsubseteq \vartheta$ and Claims \ref{ceilbound1}, \ref{ceilbound2} hold. 

Let us use $\# S$ to denote the cardinality of the set $S$ and work by induction on $\# \mathop{\rm dom}^+(r)$. The base case, where $\mathop{\rm dom}^+(r)=\{0\}$, is trivial, as $r\sqsubseteq \vartheta$ becomes $r(0)\leq\vartheta$, and clearly $\vartheta=r(0)$ satisfies the required properties.

For the inductive step, let $\eta=\min(\mathop{\rm dom}^+(r)\setminus\{0\})$ and consider $\widetilde r$ given by $\widetilde r(\xi)=r(\eta+\xi)$ whenever the latter is defined, so that $\widetilde r$ is just $r$ `shifted' by $\eta$. By induction hypothesis $\widetilde \vartheta=\lceil\widetilde r\rceil$ is defined and satisfies both claims.

To find $\vartheta=\lceil r\rceil$, define $\gamma=0$ if $0\not\in\mathop{\rm dom}(r)$ and $\gamma=r(0)+1$ otherwise, and set $\vartheta=\gamma+e^\eta\widetilde\vartheta$. Let us begin by showing that $r\sqsubseteq\vartheta$. First note that $r(0)<\vartheta$ in the case that $0\in{\rm dom}(r)$. Meanwhile, for $\xi\in(0,\lambda)\cap\mathop{\rm dom}(r)$ we have that
\[r(\xi)=\widetilde r(-\eta+\xi)<\ell^{-\eta+\xi}\widetilde\vartheta=\ell^{\xi}e^\eta\widetilde \vartheta= \ell^{\xi}(\gamma+e^\eta\widetilde \vartheta)=\ell^{\xi}\vartheta.\]
A similar argument shows that $\ell^\lambda\vartheta=r(\lambda)$ for $\lambda=\max({\rm dom}(r))$, thus establishing that $r\sqsubseteq\vartheta$. 

It remains to check that if $r\sqsubseteq \alpha$, then Claim \ref{ceilbound2} is satisfied by $\vartheta$. If $\widetilde\vartheta=0$ then this is obvious, for if $0\not\in\mathop{\rm dom}(r)$ then $\vartheta=e^\eta\widetilde\vartheta=0$, and if $0\in\mathop{\rm dom}(r)$ then $\vartheta=r(0)+1$ and therefore $\le^0\vartheta=r(0)+1\leq \alpha$, whereas for $\xi>0$ we have that $\le^\xi\vartheta=0\leq \le^\xi\alpha$. Hence we may assume $\widetilde\vartheta>0$.

Pick $\xi\leq \max(\mathop{\rm dom}(r))$. If $\xi\geq \eta$, observe that $\widetilde r\sqsubseteq \ell^\eta(\alpha)$ and thus, by our induction hypothesis,
\[\le^\xi\vartheta=\le^{-\eta+\xi}\le^\eta\vartheta= \le^{-\eta+\xi}\widetilde\vartheta\leq \ell^{-\eta+\xi}\ell^\eta\alpha=\ell^\xi\alpha.\]

If $\xi\in(0,\eta)$, then $\ell^{-\xi+\eta}\ell^\xi\alpha =\ell^\eta\alpha\geq\widetilde\vartheta$. It follows by Lemma \ref{leftinv} that $\ell^\xi\alpha \geq e^{-\xi+\eta}\widetilde\vartheta$. But then, as we are assuming $\widetilde\vartheta>0$ we obtain that \[\ell^\xi\vartheta=\ell^\xi(\gamma+e^\eta\widetilde\vartheta) =\ell^\xi  e^\eta\widetilde\vartheta = e^{-\xi+\eta}\widetilde\vartheta,\]
and thus $\ell^\xi\vartheta\leq \ell^\xi\alpha$.

Finally, we must see that $\vartheta\leq \alpha$. Since $r\sqsubseteq \alpha$, we have $\alpha=\gamma+\delta$ for some $\delta\geq 0$; but from the assumption that $\widetilde\vartheta>0$ we have that $\le^\eta\alpha>0$ and thus we must have $\delta>0$ (for $\le^\eta\gamma=0$). Now, $\ell^\eta\delta=\ell^\eta\alpha \geq \widetilde\vartheta$, so once again by Lemma \ref{leftinv} we obtain $\delta\geq e^\eta\widetilde\vartheta$ and thus $\alpha\geq \gamma+e^\eta\widetilde\vartheta=\vartheta.$

Thus we may set $\lceil r\rceil=\vartheta$ and obtain all the desired properties.
\endproof

As a variant, we may be interested in the least $\vartheta$ such that $r\sqsubset \vartheta$. We may construct it as follows: let $r'$ be equal to $r$ except that, for $\lambda=\max(\mathop{\rm dom}(r))$, we set $r'(\lambda)=r(\lambda)+1$. Then, it is straightforward to check that $\lceil r'\rceil$ is the least ordinal $\vartheta$ such that $r\sqsubset \vartheta$.

\begin{examp}\label{ExSimpSeq}
Consider the simple sequence with $r(0)=\varepsilon_0,r(\omega)=\omega^2,r(\omega+1)=2$, and undefined elsewhere. Let us compute $\vartheta=\lceil r\rceil$.

Since $\omega+1$ is the greatest element of ${\rm dom}(r)$ we must have $\le^{\omega+1}\vartheta=2$, and hence $\le^{\omega}\vartheta$ is of the form $\alpha+\omega^2$. Now, we cannot take $\alpha=0$, or else we would not have $\le^\omega\vartheta>r(\omega)=\omega^2$, and the least value of $\alpha$ we could take to obtain a strictly larger value is $\omega^2$. Thus $\le^\omega\vartheta=\omega^2+\omega^2$.

In view of Lemma \ref{leftinv}, the least value of $\vartheta$ that satisfies this is $e^\omega(\omega^2+\omega^2)$. Moreover, we already have $e^\omega(\omega^2+\omega^2)>e^\omega 1=\varepsilon_0=r(0)$. Thus we may set $\vartheta= e^\omega(\omega^2+\omega^2)= \varepsilon_{\omega^2+\omega^2}$, and we have that $\lceil r\rceil= \varepsilon_{\omega^2+\omega^2}$.
\end{examp}


\section{Ranks and $d$-maps}\label{ranks}

In this section we shall consider some fundamental concepts in the study of scattered spaces. We omit the proofs of those results which may already be found in \cite{BeklGabel:2011}.

Given a topological space $\langle X,\mathcal T\rangle$, we may iterate the corresponding derived set operator $d:\mathcal P(X)\to \mathcal P(X)$ via the following recursion:
\begin{enumerate}
\item $d^0A=A$
\item $d^{\xi+1}A=dd^\xi A$ for all $\xi\in{\sf On}$
\item $d^\lambda A=\displaystyle\bigcap_{\xi<\lambda}d^\xi A$ for $\lambda\in{\sf Lim}$.
\end{enumerate}
If $X$ is scattered and $d^\xi A\not=\varnothing$, then $d^\xi A$ contains an isolated point $x$ and thus $x\not\in dd^\xi A=d^{\xi+1} A$. In particular, $d^\xi X\supsetneq d^{\xi+1}X$ provided $d^\xi X\not=\varnothing$.

Thus if $\#\xi>\# X$ (recall that we use $\#$ to denote cardinality), $d^\xi X=\varnothing$, which means that for any point $x\in X$ there is some ordinal such that $x\not\in d^\xi X$. This motivates our following definition:

\begin{definition}[Rank]
If $\mathfrak X=\langle X,\mathcal T\rangle$ is a scattered space and $x\in X$, we define the {\em rank} of $x$, denoted $\uprho(x)$, to be the least ordinal $\alpha$ such that $x\not \in d^{\alpha+1} X$.

We define $\uprho(\mathfrak X)=\sup_{x\in X}(\uprho(x)+1)$. This is the {\em Cantor-Bendixon rank} of $\mathfrak X$.
\end{definition}

Many times it will turn out that ranks are not too difficult to compute; the following lemma gives an example of this.

\begin{lemm}\label{basicrank}
Given an ordinal $\Theta$, let $\uprho_{0}$ be the rank function on $\Theta_{0}$ and $\uprho_{1}$ the rank function on $\Theta_{1}$.

Then, for all $\xi<\Theta$, $\uprho_{0}(\xi)=\xi$ while $\uprho_{1}(\xi)=\le\xi$.
\end{lemm}

These equalities have already appeared in \cite{BeklGabel:2011}, and it is an instructive excercise to prove them directly (by induction on $\xi$). However, we shall not provide such a proof, as they are instances of the more general Corollary \ref{rankcomp} that we will give later.

\begin{examp}
Let $\mathfrak X$ be $[0,\omega^\omega]$ with the interval topology. Then, every ordinal that is either zero or a successor is isolated; $\{0\}$ is open as is $\{\xi+1\}=(\xi,\xi+2)$ for all $\xi$. Meanwhile, limit ordinals are not isolated; for example, any neighborhood of $\omega$ contains an interval $(N,\omega+1)$ and hence a point $N+1\not=\omega$. It follows that $d[0,\omega^\omega]$ is the set of limit ordinals below $\omega^\omega+1$.

Now, $\omega$ is isolated in $d[0,\omega^\omega]$ since we have removed all natural numbers, so that, for example, $(0,\omega+1)\cap d[0,\omega^\omega]=\{\omega\}$. The same situation occurs for any ordinal of the form $\gamma+\omega$. However, $\omega^2$ is {\em not} isolated in $d[0,\omega^\omega]$, as any neighborhood of $\omega^2$ contains all elements of the form $\omega\cdot N$ for $N$ large enough. More generally, no ordinal of the form $\gamma+\omega^n$ is isolated in $d[0,\omega^\omega]=\{\omega\}$ if $n\geq 2$, and $d^2[0,\omega^\omega]$ is the set of all ordinals $\xi$ below $\omega^\omega+1$ such that $\le\xi\geq 2$.

This analysis could be carried further to see that $d^n[0,\omega^\omega]$ contains exactly those elements $\xi$ with $\le\xi\geq n$. It follows that $d^\omega[0,\omega^\omega]=\{\omega^\omega\}$, and thus $\omega^\omega$ is isolated in $d^\omega[0,\omega^\omega]$, which means that it does not belong to $d^{\omega+1}[0,\omega^\omega]$ and $\uprho\omega^\omega=\omega$. We conclude that $\uprho(\mathfrak X)=\omega+1$.
\end{examp}

As the derived set operator is central to the semantics of $\glp_\Lambda$, we need to focus on those operators that preserve it. Of course, $d$ is homeomorphism-invariant, but this class of maps is too restrictive. Meanwhile, if $f$ is merely continuous and open, it is not generally the case that $df=d$. As a simple counterexample, consider the ordinals $1$ and $\omega+1$ equipped with the interval topology, and let $f:\omega+1\to 1$ be the map that is identically zero. Of course, this is the only function between the two spaces. Further, it is easily checked to be continuous and open, yet $d(\omega+1)=\{\omega\}$ while $d 1=\varnothing.$

To this end, we need to consider $d$-maps. Recall that a space is {\em discrete} if every subset is open. If $\mathfrak X$ and $ \mathfrak Y=\langle Y,\mathcal S\rangle$ are topological spaces, a map\footnote{We write $f:\mathfrak X\to\mathfrak Y$ instead of $f:X\to Y$ when the specific topologies are relevant.} $f\colon \mathfrak X\rightarrow \mathfrak Y$ is {\em pointwise discrete} if $f^{-1}(y)$ is discrete for all $y\in Y$.

\begin{definition}[$d$-map]
Given topological spaces $\mathfrak X=\langle X,\mathcal T\rangle$ and $\mathfrak Y=\langle Y,\mathcal S\rangle$, a {\em $d$-map} from $\mathfrak X$ to $\mathfrak Y$ is a function $f\colon\mathfrak X\to \mathfrak Y$ which is continuous, open, and pointwise discrete.
\end{definition}

The property of being pointwise discrete is equivalent to the apparently stronger condition that $f^{-1}A$ is discrete whenever $A$ is (see \cite{BeklGabel:2011}). With this observation one readily obtains the following:

\begin{lemm}\label{dcomp}
The composition of $d$-maps is a $d$-map.
\end{lemm}

As an important example, the rank function itself is a $d$-map, even a ``canonical'' $d$-map in a certain sense:

\begin{lemm}\label{onlyone}
Given a scattered space $\mathfrak X$, the rank function $\uprho\colon\mathfrak X\to\uprho(\mathfrak X)_{0}$ is a $d$-map. Moreover, if $f\colon\mathfrak X\to \Theta_{0}$ is a $d$-map, it follows that $f=\uprho$.
\end{lemm}

\proof See \cite[Lemma 3.3]{BeklGabel:2011}.\endproof

Thus $\Theta_{0}$ may be seen as a final object in the category of scattered spaces with Cantor-Bendixon rank at most $\Theta$ and $d$-maps as morphisms. One immediate consequence is that $d$-maps are rank-preserving. To be precise, if $\mathfrak X,\mathfrak Y$ are scattered spaces with rank-functions $\uprho_\mathfrak X,\uprho_\mathfrak Y$ and $f:\mathfrak X\to \mathfrak Y$, we say $f$ is {\em rank-preserving} if $\uprho_\mathfrak X=\uprho_\mathfrak Y f$.

\begin{lemm}\label{LemmRankPreserve}
If $\mathfrak X,\mathfrak Y$ are scattered spaces and $f\colon\mathfrak X\to\mathfrak Y$ is a $d$-map, then $f$ is rank-preserving.
\end{lemm}

\proof
Let $\uprho_\mathfrak X$ and $\uprho_\mathfrak Y$ be the respective rank functions. The maps $\uprho_\mathfrak Y f\colon \mathfrak X\to\uprho_\mathfrak Y(\mathfrak Y)$ and $\uprho_\mathfrak X \colon\mathfrak X\to \uprho_\mathfrak X(\mathfrak X)$ are both $d$-maps by Lemmas \ref{dcomp} and \ref{onlyone}, hence also by Lemma \ref{onlyone}, they must be equal.
\endproof

\begin{lemm}\label{rankopen}
Given a scattered space $\langle X,\mathcal T\rangle$ with rank function $\uprho$ and $x\in X$, we have that if $V$ is any neighborhood of $x$, then $\uprho(V\setminus\{x\})\supseteq[0,\uprho(x))$.

Moreover, there is a neighborhood $U$ of $x$ with $\uprho(U\setminus\{x\})=[0,\uprho(x))$.
\end{lemm}

\proof
By Lemma \ref{onlyone}, $\uprho$ is a $d$-map, hence continuous and open. It follows that if $V$ is any neighborhood of $x$, then $\uprho(V)$ is a neighborhood of $\uprho x$, i.e. an initial segment containing $[0,\uprho x]$. It follows that $[0,\uprho x)\subseteq \uprho(V\setminus \{x\})$.

Now, since $x$ is isolated in $d^{\uprho (x)}X$, there must be a neighborhood $U$ of $x$ such that $U\cap d^{\uprho (x)}X=\{x\}$, and hence $\uprho(U\setminus \{x\})\subseteq [0,\uprho x)$. By the previous claim, we in fact get $\uprho(U\setminus \{x\})= [0,\uprho x)$, as desired.
\endproof

The next result shows that it is particularly easy to compute the limit points of sets that are rank-determined; that is, sets such that whenever $x\in A$ and $\uprho(x)=\uprho(y)$, then $y\in A$.

\begin{lemm}\label{rhoint}
If $\mathfrak X=\langle X,\mathcal T\rangle$ is a scattered space with rank function $\uprho$ and $S$ a set of ordinals then $d\uprho^{-1}S$ is the set of all $x\in X$ such that $\uprho x>\min S$.
\end{lemm}

\proof
This is a direct consequence of Lemma \ref{rankopen}. Indeed, if $\uprho x\leq \min S$ then there is a neighborhood $U$ of $x$ such that $\uprho(U\setminus\{x\})=[0,\uprho(x))$, and hence $x\not\in d\uprho^{-1}S$; meanwhile, if $\uprho(x)>\min S$, then given a neighborhood $U$ of $x$ we have, once again by Lemma \ref{rankopen}, that $\uprho(U\setminus\{x\})\supseteq [0,\uprho x)$, hence it must contain $\min S$. We conclude that $U\setminus\{x\}$ contains a point in $\uprho^{-1}S$ different from $x$, and since $U$ is arbitrary,  $x\in d\uprho^{-1}S$.
\endproof

There is one more extension of a scattered topology that will be useful to consider. \begin{definition}
Given a topological space $\mathfrak X=\langle X,\mathcal T\rangle$, let $\dot{\mathcal T}$ be the topology generated by $\mathcal T$ and all sets of the form $d^{\xi+1}X$ such that $\xi\in\sf On$. Then, define $\dot{\mathfrak X}=\langle X,\dot{\mathcal T}\rangle$.
\end{definition}

The following claim is a modification of a result in \cite{BeklGabel:2011}:

\begin{lemm}\label{liftgen}
If $\mathfrak X,\mathfrak Y$ are scattered spaces and $f:\mathfrak X\to \mathfrak Y$ is a $d$-map then $f:\dot{\mathfrak X}\to \dot{ \mathfrak Y}$ is also a $d$-map.
\end{lemm}

\proof
Let $\mathfrak X=\langle X,\mathcal T\rangle$ and $\mathfrak Y=\langle Y,\mathcal S\rangle$. Obviously $f$ is pointwise discrete as a map from $\dot{\mathfrak X}$ to $\dot{ \mathfrak Y}$. Let us show that $f\colon \dot{\mathfrak X} \to\dot{\mathfrak Y}$ is continuous. Suppose that $V\cap d^{\zeta+1}Y$ is an $\dot{\mathcal S}$-open set. Since $f$ is rank-preserving by Lemma \ref{LemmRankPreserve}, we have that $f^{-1}d^{\zeta+1}Y=d^{\zeta+1}X$, and thus $f^{-1}(V\cap d^{\zeta+1}X)=f^{-1}(V)\cap d^{\zeta+1}X$, which is $\dot{\mathcal T}$-open.

The argument that $f\colon \dot{\mathfrak X} \to\dot{\mathfrak Y}$ is open is very similar. Let $U\cap d^{\xi+1}X$ be a $\dot{\mathcal T}$-open set. Once again, we have that $f^{-1}d^{\xi+1}Y=d^{\xi+1}X$, and thus $f(U\cap d^{\xi+1}X)=f(U)\cap d^{\xi+1}Y$, which is $\dot{\mathcal S}$-open.
\endproof

The above result will be useful in extending constructions to successor modalities. We will also need the following lemma in order to deal with limit modalities:

\begin{lemm}\label{unions}
Let $\mathfrak X=\langle X,\vec{\mathcal T}\rangle$ and $\mathfrak Y=\langle Y,\vec{\mathcal S}\rangle$ be $\lambda$-polytopologies such that both $\vec{\mathcal T}$ and $\vec{\mathcal S}$ are increasing.

If $\lambda\in{\sf Lim}$ and $f:\mathfrak X_\xi\to\mathfrak Y_\xi$ is a $d$-map for all $\xi<\lambda$ then
\[f:\Big\langle X,\bigsqcup_{\xi<\lambda}\mathcal T_\xi\Big\rangle\to\Big \langle Y,\bigsqcup_{\xi<\lambda}\mathcal S_\xi\Big\rangle\]
is a $d$-map.
\end{lemm}

\proof
Let $\mathcal T_\lambda=\bigsqcup_{\xi<\lambda}\mathcal T_\xi$ and $\mathcal S_\lambda=\bigsqcup_{\xi<\lambda}\mathcal S_\xi$. It is obvious that $f$ is pointwise discrete as $f^{-1}(y)$ is already $0$-discrete in $\mathcal T_0$ and thus also $\lambda$-discrete\footnote{We often index topological properties by the topology they refer to, i.e. {\em $0$-discrete} means {\em discrete in $\mathcal T_0$.}}.

To check that $f$ is open, suppose that $U\subseteq X$ is open in $\mathcal T_\lambda$, so that $U=\bigcup_{\lambda<\Lambda}U_\lambda$. Then, $f(U)=\bigcup_{\lambda<\Lambda}f(U_\lambda)$, which is open in $\mathcal S_\lambda$ as it is a union of open sets. Similarly, if $V\in\mathcal S_\lambda$ and $V=\bigcup_{\lambda<\Lambda}V_\lambda$ then $f^{-1}(V)=\bigcup_{\lambda<\Lambda}f^{-1}(V_\lambda)\in \mathcal T_\lambda$ and $f$ is continuous.
\endproof


\section{Icard ambiances}\label{secic}

In this section we shall discuss {\em Icard topologies}, originally introduced in \cite{Icard:2009} for ${\sf GLP}_\omega$ and generalized to arbitrary ${\sf GLP}_\Lambda$ in \cite{glpmodels}.

Let $\mathcal I_0$ be the initial segment topology, and for $0<\lambda<\Lambda$ define a topology $\mathcal I_\lambda$ on $\Theta$ by setting, for $\lambda<\Lambda$, $\mathcal I_\lambda$ to be the topology generated by sets of the form
\[(\alpha,\beta]_\xi=\cbra \vartheta:\alpha< \le^\xi\vartheta\leq\beta\cket\]
or of the form
\[[0,\beta]_\xi=\cbra \vartheta:\le^\xi\vartheta\leq\beta\cket\]
for some $\alpha<\beta\leq{{\Theta}}$ and $\xi<\lambda$. For uniformity, we may write $[0,\beta]_\xi$ as $(-1,\beta]_\xi$, and thus we may assume all intervals to be open on the left.

We will call the resulting polytopological space ${\mathfrak {Ic}}^\Theta_\Lambda$. We will denote the derived-set operator with respect to $\mathcal I_\lambda$ by $i_\lambda$ and the ordinal $\Theta$ equipped with $\mathcal I_\lambda$ by $\Theta_\lambda$; note that there is no clash in notation in the cases $\lambda=0,1$ as the Icard topologies coincide with the initial segment and interval topologies, respectively. When we need to be more specific, we will write $\mathcal I^\Theta_\Lambda$ or $i^\Theta_\Lambda$ to indicate that the underlying set is $\Theta$.

Recall that if $r$ is a simple function and $\alpha$ is an ordinal, we write $r\sqsubset\alpha$ if $r(\xi)<\le^\xi\alpha$ for all $\xi\in {\rm dom}(r)$ and $r(\lambda)\leq \le^\lambda\alpha$. If $r$ is a simple function such that $r\sqsubset \alpha$, then we can associate an Icard-neighborhood of $\alpha$ to $r$. Namely, define
\[B_r(\alpha)=\bigcap_{\xi\in{\rm dom}(r)}(r(\xi),\le^\xi\alpha]_\xi.\]
Note that if $\lambda>\max({\rm dom}(r))$, then $B_r(\alpha)$ is $\lambda$-open. This will give us a useful way to describe ``small'' neighborhoods of $\alpha$. To be precise, given a topological space $\langle X,\mathcal T\rangle$ and $x\in X$, say a family of open sets $\mathcal N$ is a {\em neighborhood base} for $x$ if $x\in U$ for all $U\in\mathcal N$ and, given any neighborhood $V$ of $x$, there is $V'\in\mathcal N$ with $V'\subseteq V$.

\begin{lemm}\label{LemmNeighBase}
If $\Theta$ is any ordinal and $\xi<\Theta$ is any ordinal such that $\le^\lambda\xi>0$, then the sets of the form $B_r(\xi)$ with $\mathop{\rm dom}(r)\subseteq \lambda$ form a neighborhood base for $\xi$.
\end{lemm}

\proof
Every neighborhood of $\xi$ contains a set of the form $V=\bigcap_{k<K}(\alpha_k,\beta_k)_{\sigma_k}$ with $\xi\in V$ and all $\sigma_k<\lambda$. Since $\le^\lambda\xi>0$, so is $\le^{\sigma_k}\xi$ for all $k$ and we may assume that all $\alpha_k$ are different from $-1$. We may also assume that all $\sigma_k$ are distinct, for if $\sigma_j=\sigma_k$ we have that $(\alpha_j,\beta_j)_{\sigma_j}\cap (\alpha_k,\beta_k)_{\sigma_k}=(\max\{\alpha_j,\alpha_k\},\min\{\beta_j,\beta_k\})$. Thus we may define $r$ by $r(\sigma_k)=\alpha_k$, with $r(\zeta)$ undefined elsewhere. Clearly, $B_r(\xi)\subseteq V$.
\endproof

It is well-known that Icard spaces are not models of ${\sf GLP}$, but as we shall see Icard {\em ambiances} are:

\begin{definition}[Icard ambiance]
An Icard ambiance is a provability ambiance $\mathfrak X$ based on an Icard space.
\end{definition}

Icard ambiances are rather nice to work with, since the topologies are all easy to describe. We will also use {\em shifted} Icard ambiances, based on $\langle\mathcal I_{1+\lambda}\rangle_{\lambda<\Lambda}$; these are important as ${\sf GLP}_1$ is incomplete for the initial segment topology.  We will denote the shifted $\Lambda$-Icard space on $\Theta$ by $\widehat{{\mathfrak {Ic}}}^\Theta_\Lambda$.

The following useful property is a slight modification of a result from \cite{glpmodels}: 

\begin{lemm}\label{onlyyou}
Given $\xi\leq{{\Theta}}$ and $\lambda<\Lambda$, there is an $\mathcal I_\lambda$-neighborhood $U$ of $\xi$ such that whenever $\xi\not=\zeta\in U$,  $\le^\lambda\zeta<\le^\lambda\xi$.
\end{lemm}

\proof
By induction on $\lambda$.

First assume that there is $\eta<\lambda$ with $\le^{\eta}\xi=0$. By induction there is an $\eta$-neighborhood $U$ of $\xi$ such that for all $\zeta\in U$ different from $\xi$, $\le^{\eta}\zeta<\le^{\eta}\xi$, which clearly imples that $U=\{\xi\}$. Since $U$ is also a $\lambda$-neighborhood of $\xi$, the result follows.

Now suppose that $\le^\eta\xi>0$ whenever $\eta<\lambda$ and consider two subcases. If $\lambda=\alpha+1$, we have by induction hypothesis that there is an $\alpha$-neighborhood $V$ of $\xi$ such that whenever $\zeta\not=\xi$ in $V$ we have $\le^\alpha\zeta< \le^\alpha\xi$. Write $\le^\alpha\xi$ as $\eta+\omega^\beta$ and consider the $\lambda$-neighborhood $U=V\cap(\eta,\le^\alpha\xi]_\alpha$. If $\beta=0$, then $U=\{\xi\}$, for any other point $\zeta\in V$ satisfies $\le^\alpha\zeta\leq \eta$ and thus does not belong to $U$. If $\beta>0$, then $\omega^\beta=e\beta$. Suppose that $\zeta\not=\xi$ belongs to $U$. We have that $\zeta=\eta+\delta$ for some $\delta \in(0,\omega^\beta]$ while from $\zeta\in V$ we obtain $\delta<\omega^\beta=e\beta$. We then have by Lemma \ref{leftinv} that $\le\delta<\beta$ and thus $\le^\lambda\zeta=\le\delta<\beta=\le^\lambda\xi.$

Finally, if $\lambda$ is a limit ordinal, use Lemma \ref{LemmLogConverge} to find $\alpha<\lambda$ and $\rho>0$ such that $\le^\alpha\xi=e^{\omega^\rho}\le^\lambda\xi$. By induction hypothesis, there is an $\alpha$-neighborhood $U$ of $\xi$ such that whenever $\zeta\not=\xi$ in $U$ we have $\le^\alpha\zeta< \le^\alpha\xi$. Then, $U$ already satisfies the desired properties; for indeed, since $\le^\alpha\zeta<\le^\alpha\xi= e^{\omega^\rho}\le^\lambda\xi$ we also have, by Lemma \ref{leftinv}, that $\le^\lambda\zeta=\le^{\alpha+\omega^\rho}\zeta= \le^{\omega^\rho}\le^\alpha\zeta< \le^\lambda\xi$.
\endproof

Constructing $d$-maps between Icard spaces will be crucial. Fortunately, hyperlogarithms already provide important examples. For simplicity, we shall henceforth write $\le^{-\xi}$ instead of $(\le^\xi)^{-1}$.

\begin{lemm}\label{lisd}
If $\Theta,\xi,\zeta$ are ordinals, then $\le^\xi:\Theta_{\xi+\zeta}\to\Theta_{\zeta}$ is a $d$-map.
\end{lemm}

\proof
That $\le^\xi$ is pointwise discrete is an immediate consequence of Lemma \ref{onlyyou}, so it remains to show that the maps are open and continuous.

Let us first consider the case when $\zeta=0$. Let $[0,\beta]_0$ be a $0$-open set and $\vartheta\in\le^{-\xi}[0,\beta]_0$. Once again use Lemma \ref{onlyyou} to find a $\xi$-neighborhood $U$ of $\vartheta$ such that for all $\eta\in U$, $\le^\xi\eta\leq \beta$. But then, $U\subseteq \le^{-\xi}[0,\beta]_0$, and since all parameters were arbitrary we conclude that $\le^\xi:\Theta_\xi\to\Theta_0$ is continuous.

To see that it is open, let $U= \bigcap_{n\leq N}(\alpha_n,\beta_n]_{\delta_n}$ be $\xi$-open and suppose that $\vartheta\in U$. We claim that $[0,\le^\xi\vartheta]\subseteq \le^\xi U$. To see this, pick $\eta\leq\le^\xi \vartheta$. Define a simple function $r$ with $r(\delta_n)=\alpha_n$ and $r(\xi)=\eta$. Then, by Lemma \ref{ceilbound}, $\le^\xi\lceil r\rceil=\eta$ while for all $n\leq N$,
\[\alpha_n< \le^{\delta_n}\lceil r\rceil\leq\le^{\delta_n}\vartheta\leq \beta_n.\]
Thus $\lceil r\rceil\in U$, so that $\eta\in\le^\xi U$. 
Since $\eta$ was arbitrary, we conclude that $[0,\le^\xi\vartheta]_0\subseteq \le^\xi U$, and thus $\le^\xi:\Theta_\xi\to \Theta_0$ is open, as claimed.

Now we must consider $\zeta>1$. To see that $\le^\xi:\Theta_{\xi+\zeta}\to\Theta_\xi$ is continuous, note that if $\delta<\zeta$, $\le^\xi\gamma\in(\alpha,\beta)_\delta$ if and only if $\le^\delta\le^\xi\gamma=\le^{\xi+\delta}\gamma\in(\alpha,\beta)$, that is,  $\le^{-\xi}(\alpha,\beta)_{\delta}=(\alpha,\beta)_{\xi+\delta}$, which is $(\xi+\zeta)$-open.

Next, let us check that it is open. Suppose that $\gamma\in U=\bigcap_{n<N}(\alpha_n,\beta)_{\delta_n}$ where $\delta_n<\delta_{n+1}<\xi+\zeta$ and suppose that $J\leq N$ is the largest index such that $\delta_n<\xi$ for all $n<J$. Consider the $\zeta$-neighborhood
\[V=[0,\le^\xi\gamma]_0\cap\bigcap_{J\leq n<N}(\alpha_n,\le^{\delta_n}\gamma]_{-\xi+\delta_n}\]
of $\le^\xi\gamma$ and choose $\eta\in V$. Define a simple function $s$ by $s(\delta_n)=\alpha_n$ for $n<J$ and $s(\xi)=\eta$, and undefined otherwise. By Lemma \ref{ceilbound} we know that $\alpha_n<\lceil s\rceil\leq\le^{\delta_n}\gamma$ for all $n<J$, while $\le^\xi\lceil s\rceil=\eta$ and for $n\geq J$, $\le^{\delta_n}\lceil s\rceil=\le^{-\xi+\delta_n}\eta\in(\alpha_n,\beta_n]$. It follows that $\lceil s\rceil\in U$ and $\le^\xi\lceil s\rceil=\eta$, so that $\eta\in\le^\xi U$, as claimed.
\endproof

An important corollary of this is the following:

\begin{cor}\label{rankcomp}
If $\uprho_\xi$ denotes the rank with respect to $\mathcal I_\xi$, then $\uprho_\xi=\le^\xi$.
\end{cor}

\proof
Immediate from Lemma \ref{lisd} with $\zeta=0$ and Lemma \ref{onlyone}.
\endproof

With this we may also obtain another useful characterization of Icard topologies from \cite{BeklGabe:2012}:

\begin{lemm}\label{dotchar}
Given ordinals $\Theta,\xi$,
\begin{enumerate}
\item If $\xi=\zeta+1$ then $\mathcal I_\xi=\dot{\mathcal I}_\zeta$,
\item if $\xi\in\sf Lim$ then $\mathcal I_\xi=\bigsqcup_{\zeta<\xi}\mathcal I_\zeta.$
\end{enumerate}
\end{lemm}

\proof
The second claim is immediate from the definitions, so we shall check only the first.

Here we note that $\mathcal I_\xi$ is obtained from $\mathcal I_\zeta$ by adding sets of the form $(\alpha,\beta]_\xi$ as opens. In view of Lemma \ref{onlyyou} we know that $[0,\beta]_\xi$ is always $\zeta$-open, so it suffices to prove that $(\alpha,\Theta)_\zeta$ is $\dot{\mathcal I}_\zeta$-open as well. But this follows from Corollary \ref{rankcomp}, as for all $\delta$ we have that $\uprho_\zeta\delta=\le^\zeta\delta$ and hence $\delta\in[\alpha,\Theta)_\zeta$ if and only if $\uprho_\zeta\delta\geq\alpha$, i.e. if $\delta\in i^{\alpha}_\zeta [0,\Theta)$. We conclude that $(\alpha,\Theta)_\zeta=i^{\alpha+1}_\zeta [0,\Theta)$, which by definition is an element of $\dot{\mathcal I}_\zeta$.
\endproof

This gives us one further result:

\begin{lemm}\label{allicard}
If $\Theta,\Xi$ are ordinals and $f:\Theta_\alpha\to\Xi_\beta$ is a $d$-map then for all $\gamma$, $f:\Theta_{\alpha+\gamma}\to\Xi_{\beta+\gamma}$ is a $d$-map.
\end{lemm}

\proof
By a simple induction on $\gamma$. For successor $\gamma$ we use Lemma \ref{dotchar} together with Lemma \ref{liftgen}; for limit $\gamma$ we use Lemma \ref{unions}.
\endproof


\section{The simple ambiance}\label{secsimp}

It will be convenient to focus on a specific ambiance to get a feel for how these may be constructed. The simple ambiance we shall present here is not entirely central to our completeness proof since $\sf GLP$ is not complete for the class of simple ambiances, but simple sets nevertheless provide the appropriate semantics for the {\em closed fragment} ${\sf GLP}^0_\Lambda$ of ${\sf GLP}_\Lambda$, where propositional variables may not occur (only $\top$).

It turns out that a set is simple if and only if it is definable by a closed formula; we will prove one implication later. With this equivalence in mind, the set of valid formulas of ${\sf L}_\Lambda$ over the class of simple ambiances is equal to the set of validities over the closed-fragment definable sets. This logic is described in \cite{joosticard,joostthesis} and extends ${\sf GLP}_\Lambda$ by the axioms for linear frames.

\begin{definition}
A set $S\subseteq\Theta$ is {\em simple} if there exist natural numbers $N,M$ and ordinals $\alpha_{nm},\beta_{nm},\sigma_{nm}$ (with $\alpha_{nm}$ possibly equal to $-1$) such that
\[S=\bigcup_{n<N}\bigcap_{m<M}(\alpha_{nm},\beta_{nm}]_{\sigma_{nm}}.\]

If all $\sigma_{nm}\leq \lambda$, we say $S$ is {\em $\lambda$-simple}.
\end{definition}

It is an easy observation that all $\lambda$-simple sets are $\lambda$-open.

\begin{lemm}\label{icardopen}
If $S,T$ are simple sets, then $\Theta\setminus S$, $S\cap T$, $S\cup T$ and $i_\lambda S$ are simple sets; further, $i_\lambda S$ is $(\lambda+1)$-open.
\end{lemm}

\proof
We focus on showing that $i_\lambda S$ is $(\lambda+1)$-simple as the other properties use standard Boolean algebra manipulations.

Note that
\[i_\lambda\bigcup_{n<N}\bigcap_{m<M}(\alpha_{nm},\beta_{nm}]_{\sigma_{nm}}=\bigcup_{n<N}i_\lambda \bigcap_{m<M}(\alpha_{nm},\beta_{nm}]_{\sigma_{nm}},\]
so our claim will be established if we prove that $i_\lambda\bigcap_{k<K}(\alpha_k,\beta_k]_{\sigma_k}$ is always $(\lambda+1)$-simple.

Thus we suppose that
\[S=\bigcap_{k\leq K}(\alpha_k,\beta_k]_{\sigma_k}.\]
Assume that $S\not=\varnothing$, since otherwise the claim is trivial given that $i_\lambda\varnothing=\varnothing$, which is $(\lambda+1)$-simple, and let $\delta\in S$. Assume also that the $\sigma_k$'s are in increasing order and let $H$ be the largest index such that $\sigma_H< \lambda$. Let $r$ be a simple function defined by $r(\sigma_k)=\alpha_k$ for all $k<K$, $r(\sigma_K)=\alpha_K+1$ and let $\alpha_\ast=\le^{\lambda}\lceil r\rceil$.

We claim that
\[i_\lambda \bigcap_{k\leq K}(\alpha_k,\beta_k]_{\sigma_k}=(\alpha_\ast,\Theta)_{\lambda}\cap\bigcap_{k\leq H}(\alpha_k,\beta_k]_{\sigma_k}.\]
Let us begin by showing that the right-hand side is contained in the left. Let $\xi\in (\alpha_\ast,\Theta)_{\lambda}\cap\bigcap_{k\leq H}(\alpha_k,\beta_k]_{\sigma_k}$ and pick any $\lambda$-neighborhood $U$ of $\xi$; in view of Lemma \ref{LemmNeighBase}, we may assume $U$ is of the form $B_t(\xi)$ for some simple function $t$ with $\mathop{\rm dom}(t)\subseteq \lambda$. Let $r'$ be a simple function which is equal to $r$ on all $\xi<\lambda$, but $r'(\lambda)=\alpha_\ast$ and $r'$ is undefined otherwise.  Then define $s=t\sqcup r'$; we claim that $\zeta=\lceil s\rceil \in B_t(\xi)\cap S$.

First note that, by Lemma \ref{ceilbound}, $\zeta\in B_t(\xi)$. Further, also using Lemma \ref{ceilbound},
\[\alpha_k< \le^{\sigma_k}\zeta \leq \le^{\sigma_k}\xi\leq\beta_k\] for all $k\leq H$,
and $\le^\lambda\zeta=\alpha_\ast$ so that for $k>H$ we see that
\[\alpha_i<\le^{\sigma_k}\lceil r\rceil=\le^{\sigma_k} \zeta\leq\le^{\sigma_k}\delta\leq \beta_i,\]
and $\zeta\in S$.

Finally, note that $\le^{\lambda}\zeta=\alpha_\ast<\le^{\lambda}\xi$, and therefore $\zeta\not=\xi$. Since $B_t(\xi)$ was arbitrary, we conclude that $\xi\in i_\lambda S$.

Now let us show that the left-hand side is contained in the right-hand side. To do this, pick $\xi\in i_\lambda S$. For each $k\leq H$, it is easy to see that $\xi\in(\alpha_k,\beta_k]_{\sigma_k}$; otherwise, $[0,\alpha_k]_{\sigma_k}\cup(\beta_k,\Theta)_{\sigma_k}$ is a $\lambda$-neighborhood of $\xi$ which does not intersect $S$. Meanwhile, if $\le^{\lambda}\xi\leq\alpha_\ast$, by Lemma \ref{onlyyou}, there is a $\lambda$-neighborhood $V$ of $\xi$ such that if $\zeta\not=\xi$ is contained in $V$, then $\le^{\lambda}\zeta<\le^{\lambda}\xi$. But by Lemma \ref{ceilbound}, if $r\sqsubset\zeta$ then $\le^{\lambda}\zeta\geq\alpha_\ast$, which means that $V\cap(S\setminus\{\xi\})=\varnothing.$
\endproof

With this we may prove the following:

\begin{theorem}
Every simple ambiance is an Icard ambiance and thus ${\sf GLP}_\Lambda$ is sound for the class of simple $\Lambda$-ambiances.
\end{theorem}

\proof
Icard polytopologies are clearly scattered and increasing, and simple sets form a provability ambiance due to Lemma \ref{icardopen}.
\endproof

We conclude by mentioning a result relating simple sets to closed formulas. The converse claim is also true, i.e. that every simple set may be defined by a closed formula, but we shall not go into details here.

\begin{lemm}\label{issimp}
Given a closed formula $\phi$ and an Icard ambiance $\mathfrak X$ on an ordinal $\Theta$ with valuation $\lb\cdot\rb$, $\lb\phi\rb$ is a simple set.
\end{lemm}

\proof
By induction on the build of $\phi$ using Lemma \ref{icardopen} and the fact that $\lb\top\rb=(-1,\Theta)_0$ is simple.
\endproof


\section{Beklemishev-Gabelaia spaces}\label{bgsec}

One key observation when constructing GLP-spaces is that the operation $\cdot^+$ is not monotone. Thus it is possible that $\mathcal T^+$ is discrete yet for some suitable refinement $\mathcal T'$ of $\mathcal T$, $(\mathcal T')^+$ is not. In this case, it will be useful to pass to such an extension. This idea is central to the completeness proof of \cite{BeklGabel:2011}.

It remains to define what refinements are `suitable'; these are given by the following definition.

\begin{definition}[rank-preserving, limit-maximal refinement]
Let $\mathcal T\subseteq\mathcal T'$ be two topologies on a set $X$. Assume $\langle X,\mathcal T\rangle$ is scattered (so that $\langle X,\mathcal T'\rangle$ is scattered as well). Let $\uprho,\uprho'$ be the respective rank functions.

Then, $\mathcal T'$ is a
\begin{enumerate}
\item {\em rank-preserving} refinement of $\mathcal T$ if $\uprho=\uprho'$;
\item {\em limit-refinement of $\mathcal T$} if it is a rank-preserving refinement and, whenever $\uprho(\xi)\not\in \sf Lim$ and $U$ is any $\mathcal T'$-neighborhood of $\xi$, there is a $\mathcal T$-neighborhood $V$ of $\xi$ such that $V\subseteq U$;
\item {\em limit-maximal refinement of $\mathcal T$} if there is no limit-refinement $\mathcal T''$ of $\mathcal T$ such that $\mathcal T'\subsetneq\mathcal T''$.
\end{enumerate}
\end{definition}

Limit-maximal refinements are very useful for constructing GLP-spaces. The following results are proven in \cite{BeklGabel:2011} and are crucial in the construction. Recall that given $\mathfrak X=\langle X,\mathcal T\rangle$, $\dot{\mathcal T}$ is the topology generated by $\mathcal T$ and all sets of the form $d^{\xi+1}X$.

\begin{lemm}\label{LemmBG}
Let $\mathfrak X,\mathfrak Y$ be scattered spaces. Then,
\begin{enumerate}
\item \label{LemmBGOne} There exists a limit-maximal refinement of $\mathfrak X$.
\item \label{LemmBGTwo} If $\mathfrak X$ is limit-maximal then $\mathfrak X^+=\dot{\mathfrak X}$.
\item  \label{LemmBGThree} If $f:\mathfrak X\to \mathfrak Y$ is a $d$-map, $\mathfrak Y$ is limit-maximal and $\mathfrak X'$ is a limit-maximal refinement of $\mathfrak X$, then $f:\mathfrak X'\to\mathfrak  Y$ is also a $d$-map.
\item  \label{LemmBGFour}
If $f:\mathfrak X\to \mathfrak Y$ is a $d$-map and $\mathfrak Y'$ is a limit-maximal refinement of $\mathfrak Y$ then there is a limit-maximal refinement $\mathfrak X'$ of $\mathfrak X$ such that $f:\mathfrak X'\to\mathfrak  Y'$ is also a $d$-map.
\end{enumerate}
\end{lemm}

\proof
These claims are all proven in \cite{BeklGabel:2011}, where Item \ref{LemmBGOne} is Lemma 4.4, Item \ref{LemmBGTwo} is Lemma 5.1, Item \ref{LemmBGThree} is Lemma 4.6 and Item \ref{LemmBGFour} Lemma 4.7.
\endproof 

With this we have all the tools we need to construct Beklemishev-Gabelaia spaces. Below, we use $\mathcal T^-_\xi$ in the sense of Definition \ref{tminus}.

\begin{definition}
Let $\Theta$ be an ordinal.

A polytopology $\langle\mathcal T_\xi\rangle_{\xi<\Lambda}$ on $\Theta$ is a {\em Beklemishev-Gabelaia space} (BG-space) if $\mathcal T_0$ is a limit-maximal refinement of the interval topology and for all $\xi$, $\mathcal T_\xi$ is a limit-maximal refinement of $\mathcal T_\xi^-$.
\end{definition}

There is a close relationship between BG and Icard spaces:

\begin{lemm}\label{bgicard}
If $\mathfrak X$ is a $\Lambda$-BG space with topologies $\vec{\mathcal T}$ then for all $\lambda<\Lambda$, $\mathcal T_\lambda$ is a rank-preserving refinement of $\mathcal I_{1+\lambda}$.
\end{lemm}

\proof
Suppose that $\mathfrak X$ is based on an ordinal $\Theta$. Let $\uprho_\lambda$ be the rank-function with respect to $\mathcal T_\lambda$; in view of Theorem \ref{rankcomp}, the rank with respect to $\mathcal I_{1+\lambda}$ is $\le^{1+\lambda}$. We proceed to prove that $\uprho_\lambda=\le^{1+\lambda}$ by induction on $\lambda$, with the base case $\lambda=0$ being immediate from the definitions and Lemma \ref{basicrank}. For $\lambda=\xi+1$, we use Lemma \ref{LemmBG}.\ref{LemmBGTwo} to see that $\mathfrak X_\xi^+=\dot{\mathfrak X}_\xi$, and hence $\mathcal T_\lambda$ is rank-preserving over $\dot{\mathcal T}_\xi$, so we may compute ranks over $\dot{\mathcal T}_\xi$ instead of $\mathcal T_\lambda$.

By induction hypothesis, $\mathcal T_\xi$ is a rank-preserving refinement of $\mathcal I_{1+\xi}$. Since we also have $\mathcal I_{1+\lambda}=\dot{\mathcal I}_{1+\xi}$, it readily follows that $\mathcal T_\lambda$ is a refinement of $\mathcal I_{1+\lambda}$. We use a second induction on $\vartheta$ to show that $\uprho_\lambda\vartheta=\le^{1+\lambda}\vartheta$ for $\vartheta< \Theta$; that is, assume that if $\vartheta'<\vartheta$ then $\uprho_\lambda\vartheta'=\le^{1+\lambda}\vartheta'$. Use Lemma \ref{rankopen} to find a $\dot{\mathcal T}_\lambda$-neighborhood $U\subseteq [0,\vartheta]$ of $\vartheta$ with $\uprho_\lambda(U)=[0,\uprho_\lambda\vartheta)$, so that $U=V\cap (\alpha,\le^{1+\xi}\vartheta]_{1+\xi}$ for some $V\in\mathcal T_\xi$ and $\alpha<\Theta$. Let $\delta<\le^{1+\lambda}\vartheta$. Then, there is $\gamma\in (\alpha,\le^{1+\xi}\vartheta]$ with $\le\gamma=\delta$, since $\le$ maps intervals to initial segments.

But $\gamma<\le^{1+\xi}\vartheta$ and since by induction hypothesis $\mathcal T_\xi$ is rank-preserving over $\mathcal I_{1+\xi}$, there is $\eta\in V$ with $\uprho_\xi\eta=\le^{1+\xi}\eta=\gamma$. It follows that $\eta\in U$ and, by induction on $\eta<\vartheta$, $\uprho_\lambda\eta=\le^{1+\lambda}\eta$. Since $\uprho_\lambda(U)=[0,\uprho_\lambda\vartheta)$ we have that $\uprho_\lambda\vartheta>\uprho_\lambda\eta\stackrel{\rm IH}=\delta$. Since $\delta<\le^{1+\lambda}\vartheta$ was arbitrary, we conclude that $\uprho_\lambda\vartheta\geq \le^{1+\lambda}\vartheta$. The inequality $\uprho_\lambda\vartheta\leq \le^{1+\lambda}\vartheta$ follows from the fact that $\mathcal T_\lambda$ refines $\mathcal I_{1+\lambda}$, and hence the two are equal.

If $\lambda$ is a limit ordinal, first note that $1+\lambda=\lambda$, which will simplify some expressions. We have that $\mathcal I_\lambda=\bigsqcup_{\xi<\lambda}\mathcal I_\xi$ whereas $\mathcal T_\lambda\supseteq \mathcal T^-_\lambda=\bigsqcup_{\xi<\lambda}\mathcal T_\xi$, so by induction $\mathcal T_\lambda$ is a refinement of $\mathcal I_\lambda$. It remains to show that it is rank-preserving.

Observe that $\mathcal T_\lambda$ is rank-preserving over $\mathcal T^-_{\lambda}=\bigsqcup_{\xi<\lambda} \mathcal T_\xi$, so it suffices to compute ranks over $\mathcal T^-_{\lambda}$. Pick any basic $\mathcal T^-_{\lambda}$-neighborhood $U$ of $\vartheta$, so that $U\in \mathcal T_\xi$ for some $\xi<\lambda$, and $\delta<\le^{\lambda}\vartheta$. We may assume $U\subseteq [0,\vartheta]$. By induction on $\xi<\lambda$, $\uprho_\xi U\supseteq[0,\le^{1+\xi}\vartheta]$. From $\delta<\le^\lambda\vartheta= \le^{-(1+\xi)+\lambda}\le^{1+\xi}\vartheta$ and Lemma \ref{leftinv} we obtain $e^{-(1+\xi)+\lambda}\delta<\le^{1+\xi}\vartheta$ and hence there is $\gamma\in U$ with $\uprho_\xi\gamma= e^{-(1+\xi)+\lambda}\delta$. But then, $\uprho_\lambda\gamma\in\uprho_\lambda U$, and by induction on $\gamma<\vartheta$ we have
\[\uprho_\lambda\gamma=\le^{\lambda} \gamma=\le^{-(1+\xi)+\lambda} \le^{1+\xi}\gamma= \le^{-(1+\xi)+\lambda}e^{-(1+\xi)+\lambda}\delta=\delta.\]
Since $U$ was arbitrary it follows that $\uprho_\lambda\vartheta\geq\le^{1+\lambda}\vartheta$, and hence the two are equal.
\endproof

Thus the analogue of Theorem \ref{rankcomp} also holds for BG-spaces, although here we obtain $\uprho_\xi=\le^{1+\xi}$. In fact, the above result motivates our focusing on rank-preserving extensions of Icard spaces. Both BG-spaces and shifted Icard ambiances are examples of regular polytopologies, in the sense of the following definition:

\begin{definition}[regular space]
A $\Lambda$-space $\mathfrak X$ with topologies $\vec{\mathcal T}$ is {\em regular} if for all $\lambda<\Lambda$, $\mathcal T_\lambda$ is a limit-refinement of $\mathcal I_{1+\lambda}$.
\end{definition}

It remains to show that BG-spaces actually {\em exist}. This can be done via a non-constructive proof:

\begin{lemm}\label{isbg}
Given ordinals $\Theta,\Lambda$, there exists a BG-space $\langle \Theta,\langle\mathcal T_\lambda\rangle_{\lambda<\Lambda}\rangle$.
\end{lemm}

\proof
By a straightforward induction on $\Lambda$ using Lemma \ref{LemmBG}.\ref{LemmBGOne} to find a limit-maximal refinement of $\mathcal T^-_\Lambda$; we remark that such a refinement is found using Zorn's lemma and hence the resulting space is not given constructively. The case for $\Lambda = \omega$ was first proven in \cite{BeklGabel:2011}.
\endproof

Below, a polytopology $\langle X,\langle\mathcal T_\lambda\rangle_{\lambda<\Lambda}\rangle$ is {\em based} on a topological space $\langle X,\mathcal T\rangle$ if $\mathcal T_0$ is a limit-extension of $\mathcal T$. 

\begin{lemm}\label{polylift}
Suppose that $\Theta$ is an ordinal, $\mathfrak Y$ a $\Lambda$-BG-space and $f: \Theta_1\to \mathfrak Y_0$ a $d$-map.

Then, there exists a BG-space $\mathfrak X$ based on $\Theta_1$ such that $f:\mathfrak X_\lambda\to\mathfrak Y_\lambda$ is a $d$-map for all $\lambda<\Lambda$.
\end{lemm}

\proof
We proceed by induction on $\lambda$, assuming we have constructed $\mathcal T_\xi$ for $\xi<\lambda$. If $\lambda=0$, $f:\Theta_1\to\mathfrak Y_0$ is a $d$-map by assumption so by Lemma \ref{LemmBG}.\ref{LemmBGFour} there is a limit-maximal extension $\mathfrak X_0$ of $\Theta_1$ such that $f:\mathfrak X_0\to \mathfrak Y_0$ is a $d$-map.

If $\lambda=\xi+1$ then by Lemma \ref{LemmBG}.\ref{LemmBGTwo}, $\mathfrak X^+_\xi=\dot{\mathfrak X}_\xi$ and $\mathfrak Y^+_\xi=\dot{\mathfrak Y}_\xi$ so that by Lemma \ref{liftgen}, $f:{\mathfrak X}^+_{\xi}\to{\mathfrak Y}^+_\xi$ is a $d$-map. Then, once again by Lemma \ref{LemmBG}.\ref{LemmBGFour}, we can extend $\mathfrak X_{\xi}^+$ to a limit-maximal space $\mathfrak X_\lambda$ making $f:\mathfrak X_\lambda\to\mathfrak Y_\lambda$ a $d$-map.

Finally, if $\lambda\in{\sf Lim}$, we proceed as above, using Lemma \ref{unions}.
\endproof

BG-spaces and Icard ambiances can sometimes be united into a single structure. We call these {\em idyllic ambiances:}

\begin{definition}[idyllic ambiance]
A shifted Icard $\Lambda$-ambiance $\mathfrak X=\langle\Theta,\vec{\mathcal T},\mathcal A\rangle$ is {\em idyllic} if there is a BG polytopology on $\mathfrak X$ with derived set operators $d_\lambda$ such that, for all $\lambda<\Lambda$, $d_\lambda\upharpoonright\mathcal A=i_{1+\lambda}\upharpoonright \mathcal A$.
\end{definition}

The purpose of these ambiances is to ``kill two birds with one stone'', since any model based on an idyllic ambiance may be regarded both as an Icard model and a BG model. It will also be curious to observe that in our completeness proof, we shall construct BG-models and Icard models {\em with the same valuations.}


\section{Reductive maps}\label{secred}

A fundamental technique in building models of $\sf GLP$ consists of ``pulling back'' valuations from a previously constructed model using well-behaved maps.  If $\mathfrak X,\mathfrak Y$ are scattered polytopologies, $\lb\cdot\rb_\mathfrak Y$ is a valuation on $\mathfrak Y$ and $f\colon \mathfrak X\to\mathfrak Y$, then we may define a new valuation $\lb\cdot\rb_\mathfrak X$ by setting $\lb p\rb_\mathfrak X=f^{-1}\lb p\rb_\mathfrak Y$. Specifically, we want $\mathfrak Y$ to be of the form $\Xi_1$ for some ordinal $\Xi$, since then we can borrow from the completeness of ${\sf GL}$ for the class of ordinals with the interval topology (see \cite{glcomplete}).

This idea is used in \cite{BeklGabel:2011} to prove the topological completeness of ${\sf GLP}_\omega$ using $f=\le$. One can generalize this to arbitrary ${\sf GLP}_\Lambda$ using hyperlogarithms, but it does not lead to optimal results; because of this, in this section we broaden our arsenal of useful functions by introducing {\em reductive maps.}

One key property of hyperlogarithms is that they are determined by the lower hyperlogarithms; if $\Lambda$ is additively indecomposable and $\lambda<\Lambda$, then $\le^\Lambda\xi=\le^\Lambda\zeta$ whenever $\le^\lambda\xi=\le^\lambda\zeta$. To be precise, for functions $f\colon X\to Y$ and $g\colon X\to Z$, say $f$ is {\em $g$-determined} if $f(x)=f(y)$ whenever $g(x)=g(y)$; then, it is straightforward to check that $\le^\Lambda$ is $\le^\lambda$-determined when $\lambda<\Lambda$. This property will be extremely useful later; in fact, even a weaker, {\em local} version will already turn out to be quite powerful. For this, if $f,g$ are as above and $\mathcal T$ is a topology on $X$, say $f$ is {\em locally $g$-determined} if for every $x\in X$ there is a neighborhood $U$ of $x$ such that $f\upharpoonright U$ is $g$-determined. If $g=\le^\lambda$ we will say simply {\em $\lambda$-determined} instead of {\em $\le^\lambda$-determined.}

\begin{definition}
Suppose that $\Theta,\Xi$ are ordinals. We say a $d$-map $f\colon \Theta_{1+\Lambda}\to \Xi_1$ is {\em $\Lambda$-reductive} if for all $\lambda<1+\Lambda$, $f$ is $\mathcal I_\lambda$-locally $\lambda$-determined. 
\end{definition}

In other words, given $\vartheta<\Theta$ and $\lambda<1+\Lambda$, there is a $\lambda$-neighborhood $U$ of $\vartheta$ such that if $\xi,\zeta\in U$ and $\le^\lambda\xi=\le^\lambda\zeta$, then $f\xi=f\zeta$. As mentioned above, reductive maps generalize hyperlogarithms:

\begin{lemm}\label{LemmLogRed}
Given ordinals $\Theta$ and $ \Lambda>0$, the hyperlogarithm \[\le^{\Lambda}:(e^{\Lambda+1}\Theta+1)_{\Lambda+1}\to(\Theta+1)_1\]
is $(-1+\Lambda+1)$-reductive.
\end{lemm}

\proof
In view of Lemma \ref{lisd}, $\le^\Lambda$ is a $d$-map. Moreover, if $\lambda<\Lambda+1$ we have that $\le^{\Lambda}=\le^{-\lambda+\Lambda}\le^{ \lambda}$, so indeed given $\xi$ we have that $\le^\Lambda$ is ${\lambda}$-determined on all of $[0,\Theta]$, which is clearly $\lambda$-open.
\endproof

Thus we have reductive maps when $\Lambda$ is a successor, but for limit $\Lambda$ we will need to look elsewhere. Later in this section we will construct $d$-maps between $(e^\Lambda\Theta+1)_\Lambda$ and $(\Theta+1)_1$ when $\Lambda\in\sf Lim$, but first let us discuss some of the properties of reductive maps.

Reductive maps will be particularly important in the study of idyllic ambiances. Note that on such ambiances, $d_\lambda\upharpoonright\mathcal A=i_{1+\lambda}\upharpoonright \mathcal A$ must hold only for a {\em specific} BG-topology. But there are sets $S$ such that $d_\lambda S=i_{1+\lambda} S$ whenever $d_\lambda$ is based on {\em any} BG-polytopology. We will say such a set $S$ is {\em $\lambda$-absolute}.

\begin{definition}
Let $\Theta$ be an ordinal. A set $A\subseteq\Theta$ is {\em $\lambda$-absolute} if for every BG $(\lambda+1)$-polytopology $\vec{\mathcal T}$ on $\Theta$ we have that $d_\lambda A=i_{1+\lambda}A$.
\end{definition}

A very easy example of a $\lambda$-absolute set is the empty set, since $d\varnothing=\varnothing$ no matter what topology $d$ is defined by. Lemma \ref{rhoint} gives us a way to construct more interesting $\lambda$-absolute sets, for given any set $S$ we know that $\le^{-(1+\lambda)}S$ is $\lambda$-absolute. More generally, any reductive map gives rise to absolute sets:

\begin{lemm}\label{absolutely}
If $f$ is $\Lambda$-reductive, $\lambda<\Lambda$ and $A$ is any set, then $f^{-1}A$ is $\lambda$-absolute.
\end{lemm}

\proof
Let $f\colon\Theta\to\Xi$ be $\Lambda$-reductive, let $\lambda<\Lambda$ and let $\langle d_\xi\rangle_{\xi<\Lambda}$ be the derived-set operators for a BG polytopology $\vec{\mathcal T}$ on $\Theta$. Note that the rank function for $\mathcal T_\lambda$ is given by $\le^{1+\lambda}$.

Let $\vartheta<\Theta$ and pick a $\mathcal I_{1+\lambda}$-neighborhood $U$ of $\vartheta$ such that $f$ is $(1+\lambda)$-determined on $U$.

Meanwhile, letting $B=U\cap f^{-1}(A)$ we claim that
\[B=U\cap \le^{-(1+\lambda)}\le^{1+\lambda}B.\]
To see the left-to-right inclusion, observe that $B\subseteq \le^{-(1+\lambda)}\le^{1+\lambda}B$ and, since $B=U\cap f^{-1}(A)$, it follows that $B\subseteq U\cap \le^{-(1+\lambda)}\le^{1+\lambda}B$.

For the other inclusion, suppose that $\xi\in U\cap \le^{-(1+\lambda)}\le^{1+\lambda}B$. Then, $\le^{1+\lambda}\xi\in \le^{1+\lambda}B$, that is, there is $\xi'\in B$ such that $\le^{1+\lambda}\xi=\le^{1+\lambda}\xi'$. But since we have that $\xi,\xi'\in U$ and $f$ is $(1+\lambda)$-determined on $U$, it follows that $f\xi=f\xi'$. From $\xi'\in B$ we obtain $f\xi'\in A$ and thus $f\xi\in A$, i.e. $\xi\in U\cap f^{-1}A=B$, as desired.

Moreover, since $U$ is $\mathcal I_{1+\lambda}$-open, we have that it is $\mathcal T_\lambda$-open as well and thus
\[d_\lambda B=U\cap d_\lambda \le^{-(1+\lambda)}\le^{1+\lambda}B,\]
and similarly
\[i_{1+\lambda}B=U\cap i_{1+\lambda} \le^{-(1+\lambda)}\le^{1+\lambda}B.\]

Since $\le^{1+\lambda}$ is the rank both on $\mathcal T_\lambda$ and $\mathcal I_{1+\lambda}$, by Lemma \ref{rhoint} we see that
\begin{align*}
 d_\lambda \le^{-(1+\lambda)}\le^{1+\lambda}B&=i_{1+\lambda} \le^{-(1+\lambda)}\le^{1+\lambda}B\\
 &=\{\vartheta\leq\Theta:\le^{1+\lambda}\vartheta>\min \le^{1+\lambda}B\},
\end{align*}
from which we obtain that $\vartheta\in d_\lambda (U\cap f^{-1}(A))$ if and only if $\le^{1+\lambda}\vartheta>\min\le^{1+\lambda}B$ if and only if $\vartheta \in i_{1+\lambda}(U\cap f^{-1}(A))$. Since $\vartheta$ was arbitrary, $d_\lambda f^{-1}A=i_{1+\lambda}f^{-1}A$, as claimed.
\endproof

Another nice property of reductive maps is that they behave well with respect to extensions of limit topologies:

\begin{lemm}\label{stilld}
Suppose that that $\Theta,\Xi$ and $\Lambda$ are ordinals with $\Lambda\in\sf Lim$, $f:(\Theta+1)_\Lambda\to(\Xi+1)_1$ is a $\Lambda$-reductive map and $\langle \mathcal T_\lambda\rangle_{\lambda\leq\Lambda}$ is a regular polytopology on $\Theta+1$ with $\mathcal T_\Lambda=\bigsqcup_{\lambda<\Lambda}\mathcal T_\lambda$. Let $\mathfrak X=\langle\Theta+1,\vec{\mathcal T}\rangle$ be the resulting $(\Lambda+1)$-space.

Then, $f:\mathfrak X_\Lambda\to\Xi_1$ is a $d$-map.
\end{lemm}

\proof
Clearly $f$ is continuous and pointwise discrete, so let us check that it is open. Suppose that $U$ is a $\mathcal T_\Lambda$-neighborhood of a point $\vartheta<\Theta$, so that it is a $\mathcal T_\lambda$-neighborhood of $\vartheta$ for some $\lambda<\Lambda$. Pick an $\mathcal I_\Lambda$-neighborhood $D$ of $\vartheta$ such that $f$ is $\le^{1+\lambda}$-determined on $D$ and $\le^{1+\lambda}(D\setminus\{\vartheta\})=[0,\le^{1+\lambda}\vartheta)$; the first condition can be met because $f$ is $\Lambda$-reductive, the second by Lemmas \ref{rankopen} and \ref{rankcomp}.

We claim that $f(U\cap D)=f(D)$. Indeed, if $\zeta\in f( D)$, then $\zeta=f(\delta)$ for some $\delta\in D$. But since $\uprho_\lambda\vartheta=\le^{1+\lambda}\vartheta\geq\uprho_\lambda\delta$, there is some $\delta'\in U\cap D$ with \[\le^{1+\lambda}\delta' =\uprho_\lambda\delta'=\uprho_\lambda\delta=\le^{1+\lambda}\delta,\]
and hence $f(\delta')=f(\delta)=\zeta$. Since $\zeta$ was arbitrary, the claim follows.

Now, $f$ is a $d$-map with respect to $\mathcal I_\Lambda$, so that $f(D)$ (and hence $f(U\cap D)$) is open, as desired.
\endproof

Not all reductive maps are given by hyperlogarithms. Let us now construct another interesting example. Here we will work with {\em fundamental sequences;} that is, we assume that to each (small enough) countable limit ordinal $\xi$ we have assigned a sequence of ordinals $\langle\xi[n]\rangle_{n<\omega}$ with the property that $\xi[n]<\xi[n+1]$ for all $n$ and $\xi=\lim_{n\to \omega}\xi[n]$. If $\xi=\zeta+1$, we will define $\zeta=\xi[n]$ for all $n$.

Suppose that $\Lambda$ is infinite and additively indecomposable. If $\vartheta<e^\Lambda\Theta$ is any ordinal that is not in the range of $e^\Lambda$, there exists a value of $N$ such that $\le^{\Lambda[N]}\vartheta \leq e^\Lambda(\Theta[N])$; if $\Theta$ is a successor ordinal this is essentially Lemma \ref{falls}, otherwise $e^\Lambda\Theta=\lim_{n\to\omega}e^\Lambda(\Theta[n])$ (since $e^\Lambda$ is normal), and hence for some value of $N$ we already have that $\vartheta<e^\Lambda(\Theta[n])$.
We will denote the smallest such value of $N$ by $N^\Theta_\Lambda(\vartheta)$.

\begin{definition}
Given countable ordinals $\Theta,\Lambda$ with fundamental sequences $\langle \Theta[n]\rangle_{n<\omega}$, $\langle \Lambda[n]\rangle_{n<\omega}$ and $\vartheta<e^\Lambda\Theta$, we define $N=N^\Theta_\Lambda(\vartheta)$ to be the least natural number $N$ such that $\le^{\Lambda[N]}\vartheta \leq e^\Lambda(\Theta[N])$.
\end{definition}

\begin{examp}\label{ExNOne}
Suppose that $\Theta=2$, so that $\Theta[n]=1$ for all $n<\omega$, and $\Lambda=\omega$. A standard fundamental sequence we may take for $\Lambda$ is $\Lambda[n]=n$. Then, $e^\Lambda\Theta=\varepsilon_1$, while $e^\Lambda(\Theta[n])=\varepsilon_0$ for all $n<\omega$. Let $\vartheta=\omega^{\varepsilon_0+1}=e(e^\omega 1+1)$. Then, $\le^0\vartheta=\omega^{\varepsilon_0+1}>\varepsilon_0=e^\Lambda\Theta[0]$, $\le^1\vartheta=\varepsilon_0+1>e^\Lambda\Theta[1]$, but $\le^2\vartheta=\le\le\vartheta=\le(\varepsilon_0+1)=0$.  Thus, $N^2_\omega(\omega^{\varepsilon_0+1})=2$.

Meanwhile, if instead $\vartheta=\varepsilon_0$ we already have $\le^0\varepsilon_0=\varepsilon_0\leq e^\Lambda(\Theta[0])$, so $N^2_\omega(\varepsilon_0)=0$.
\end{examp}

\begin{examp}\label{ExNTwo}
Now suppose that $\Theta=\Lambda=\omega$, again with fundamental sequence $\omega[n]=n$. Suppose further that $\vartheta=\omega^{\varepsilon_2\cdot 3}=e((e^\omega 3)\cdot 3)$.

Then, $\le^{0}\vartheta = \omega^{\varepsilon_2\cdot 3} > 0= e^\omega(\Theta[0])$, $\le^{1}\vartheta = \varepsilon_2\cdot 3 > \varepsilon_0 = e^\Lambda(\Theta[1])$, and $\le^{2}\vartheta = \varepsilon_2 > \varepsilon_1 = e^\Lambda(\Theta[2])$. After this, $\le^{n}\vartheta = \varepsilon_2$ for all $n\geq 2$. However, $\Theta[n]$ is still increasing, so that $e^\Lambda(\Theta[3])=\varepsilon_2\geq \le^{3}\vartheta$, and thus $N^\omega_\omega(\omega^{\varepsilon_2\cdot 3})=3$.
\end{examp}

We will use the parameter $N^\Theta_\Lambda(\vartheta)$ to ``diagonalize'' and in this way define a new family of reductive maps. This parameter will be used to partition $[0,e^\Lambda]$ into countably many sets.

\begin{lemm}\label{deltalemm}
If $\Theta,\Lambda$ are ordinals such that $\Lambda$ is infinite and additively indecomposable then for every $N>1$, the set
\[\Delta^\Theta_\Lambda[N]=\{\vartheta<e^\Lambda\Theta: N^\Theta_\Lambda(\vartheta)=N\}\]
is $(\Lambda[N]+1)$-simple and $\Lambda[N]$-open.

Further,
\begin{equation}\label{secondclaim}
\ell^{\Lambda[N]} \Delta^\Theta_\Lambda[N] =[0,e^\Lambda(\Theta[N])].
\end{equation}
\end{lemm}

\proof
Let $\alpha=e^\Lambda(\Theta[N])$ and $\beta=e^\Lambda(\Theta)$.

We have that
\[\Delta^\Theta_\Lambda[N]=\bigcap_{n<N}(\alpha,\beta)_{\Lambda[n]}\cap[0,\alpha]_{\Lambda[N]},\]
which is $(\Lambda[N]+1)$-simple and, in view of Lemma \ref{onlyyou}, it is $\Lambda[N]$-open.

Now, let $\delta\leq \alpha$ and consider the simple function given by $r(\lambda[n])= \alpha$ for $n<N$ and $r(\lambda[N+1])=\delta$. By Lemma \ref{ceilbound}, $\lceil r\rceil\in\Delta^\Theta_\Lambda[N]$ and $\le^{\Lambda[N]}\lceil r\rceil=\delta.$ Since $\delta\leq \alpha=e^\Lambda(\Theta[N])$ was arbitrary, we conclude that (\ref{secondclaim}) holds.
\endproof

Recall that a family of sets $\mathcal U$ forms a {\em neighborhood base} for $x$ if every element of $\mathcal U$ is a neighborhood of $x$ and for every neighborhood $V$ of $x$ there exists $U\in\mathcal U$ such that $U\subseteq V$.

\begin{lemm}
Let $\Lambda$ be countable and additively indecomposable and $\Theta$ be any ordinal.

Then, the sets
\[D^\Theta_\Lambda[N]= (e^\Lambda(\Theta[N]),e^\Lambda(\Theta)\big]_{{\Lambda[N]}}\]
form a $\Lambda$-neighborhood base for $e^\Lambda\Theta$.

Further, $D^\Theta_\Lambda[N]=\{e^\Lambda\Theta\}\cup\bigcup_{n>N}\Delta^\Theta_\Lambda[n]$.
\end{lemm}

\proof
Clearly $e^\Lambda\Theta\in D^\Theta_\Lambda[N]$ for all $N$ and $D^\Theta_\Lambda[N]$ is open.

Now, let $B_r(e^\Lambda\Theta)$ be any basic $\Lambda$-neighborhood of $e^\Lambda\Theta$; we need to find a smaller neighborhood of the form $D^\Theta_\Lambda[N]$. Let us consider two cases.

First assume $\Theta=\Theta'+1$, so that $\Theta[n]=\Theta'$ for all $n$.

We have by Proposition \ref{prophyp}.\ref{prophypfive} that
\[e^{\Lambda}\Theta= \lim_{n\to\omega}e^{\Lambda[n]}(e^{\Lambda}\Theta'+1)\]
and hence for some value of $N$ we have that $e^{\Lambda[N]}(e^{\Lambda}\Theta'+1)>r(\xi)$ for all $\xi$.

Let $M>N$ be large enough so that $\Lambda[M]>\xi+\Lambda[N]$ for all $\xi\in{\rm dom}(r)$. Then, we claim that $D^\Theta_\Lambda[M]\subseteq B_r(\vartheta)$; for indeed, if $\xi\in {\rm dom}(r)$ and $\delta\in D^\Theta_\Lambda[M]$ then
\[\le^{\Lambda[M]}\delta=\le^{-\xi+\Lambda[M]} \le^\xi\delta\geq e^\Lambda(\Theta')+1\]
so that by Lemma \ref{leftinv}
\[\le^\xi\delta\geq e^{-\xi+\Lambda[M]}(e^{\Lambda}(\Theta')+1)\geq e^{\Lambda[N]}(e^{\Lambda}(\Theta')+1)\geq r(\xi).\]

Thus, $D^\Theta_\Lambda[M]\subseteq B_r(e^\Lambda\Theta)$.

If $\Theta$ is a limit ordinal, the argument is somewhat simpler; pick $M$ so that $\Theta[M]>e^\xi r(\xi)$ for all $\xi\in\mathop{\rm dom}(\xi)$. Then, it is easy to check that $D^\Theta_\Lambda[M]\subseteq B_r(e^\Lambda\Theta)$.

Finally, to see that $D^\Theta_\Lambda[N]=\{e^\Lambda\Theta\}\cup\bigcup_{n>N}\Delta^\Theta_\Lambda[n]$, note that every $  \xi< e^\Lambda\Theta$ lies in $\Delta^\Theta_\Lambda[n]$ for some $n$, and $n>N$ if and only if $\xi\in D^\Theta_\Lambda[N]$.
\endproof

On $[0,e\Theta]$ we shall consider a different class of neighborhoods. Define
\[\sigma_\Theta[N]=\sum_{n<N}(e(\Theta[n])+1)\]
and $\Sigma_\Theta[N]=[\sigma_\Theta[N],\sigma_\Theta[N+1])$.

Similarly, define $S_\Theta[N]=[\sigma_\Theta[N],e\Theta]$.

Then we have that:

\begin{lemm}
The sets $\{\Sigma_\Theta[n]:{n<\omega}\}$ form a partition of $(0,e\Theta)$ into $1$-open sets.

Further, the sets $\{S_\Theta[n]:n<\omega\}$ form a $1$-neighborhood base for $e\Theta$, and
\[S_\Theta[N]=\{e(\Theta)\}\cup\bigcup_{n>N}\Sigma_\Theta[n].\]
\end{lemm}

\proof
Note that despite its formal appearance the set $\Sigma_\Theta[N]$ is always open since $\sigma_\Theta[N]$ is always a successor ordinal or zero. The rest of the claims are obvious if we show that $\langle\sigma_\Theta[n]\rangle_{n<\omega}$ is unbounded in $e\Theta$.

Here we consider two cases; if $\Theta=\Theta'+1$, then $\sigma_\Theta[n]=e(\Theta')n$ for all $n<\omega$, and
\[e(\Theta'+1)=\omega^{\Theta'+1}=\lim_{n\to\omega} (\omega^{\Theta'}n+1)= \lim_{n\to\omega}(e(\Theta')+1)n.\]
But $(e(\Theta')+1)N=\sigma_\Theta[N]$.

Meanwhile, if $\Theta\in\sf Lim$, then $\sigma_\Theta[N+1]=e(\Theta[N])+1$ (all previous terms cancel) and $e\Theta=\lim_{n\to\omega}e(\Theta[N])+1$.\endproof

With this we may define the following maps:

\begin{definition}
Given countable ordinals $\Lambda,\Theta$ such that $\Lambda$ is infinite and additively indecomposable, we will define a function \[\reductive_\Lambda^\Theta:(e^{\Lambda}(\Theta)+1)\to (e(\Theta)+1)\]
assuming $\reductive^{\Theta'}_{\Lambda}$ is defined whenever $\Theta'<\Theta$.

First define $\reductive^\Theta_\Lambda e^{\Lambda}\Theta=e\Theta$.

Then, for $\xi<e^{\Lambda}\Theta$, set $N=N_\Lambda^\Theta(\xi)$ and
\[\reductive^\Theta_\Lambda\xi=\sigma_\Theta[N] +\reductive^{\Theta[N]}_ {\Lambda}\le^{\Lambda[N]} \xi.\]
\end{definition}

To illustrate the above recursion, let us show that $\reductive^\Theta_\Lambda 0$ is always zero.

\begin{lemm}\label{LemmZero}
Given arbitrary $\Theta,\Lambda$, $\reductive^\Theta_\Lambda 0=0$.
\end{lemm}

\proof
Observe that $N^\Theta_\Lambda ( 0)$ is always zero since $0\leq e^\Lambda(\Theta[0])$ independently of $\Theta$ or $\Lambda$, while $\sigma_\Theta[0]$ is always zero as well since it is an empty sum. With this, we may proceed by induction on $\Theta$; for the base case, we see that $e^\Lambda 0=0$ so $\reductive^0_\Lambda 0=e 0$. For the inductive step we have that
\[
\reductive^{\Theta}_\Lambda 0 =\sigma_\Theta[0]+\reductive^{\Theta[0]}_\Lambda 0\stackrel{\rm IH}=0+0,
\]
where we are using our induction hypothesis on $\Theta[0]<\Theta$. Thus $\reductive^{\Theta}_\Lambda 0=0$ for all $\Theta,\Lambda$, as claimed.
\endproof

\begin{examp}
In Example \ref{ExNOne} we set $\Theta=2$, $\Lambda=\omega$ and showed that $N^2_\omega(\vartheta)=2$, where $\vartheta=\omega^{\varepsilon_0+1}$. Let us use this to compute $\reductive^2_\omega\omega^{\varepsilon_0+1}$. Observe that $e(\Theta)=\omega^2$, while $e^{\Lambda}\Theta=e^\omega 2=\varepsilon_1$, so $\reductive^2_\omega\colon [0,\varepsilon_1]\rightarrow [0,\omega^2]$.

First we must set $N=N^\Theta_\Lambda(\vartheta)=2$. Since $\Theta[n]$ is the constant $1$, we have that \[\sigma_2[2]=e(2[0])+1+e(2[1])+1=\omega\cdot 2+1.\] 
Meanwhile, $\Lambda[2]=\omega[2]=2$, so $\le^{\Lambda[N]}\vartheta=\le^2\omega^{\varepsilon_0+1}= \le({\varepsilon_0+1})=0$. Thus,
\[\reductive^\Theta_\Lambda\vartheta=\sigma_\Theta[N]+ \reductive^{\Theta[N]}_ {\Lambda}\le^{\Lambda[N]}\vartheta =\omega\cdot 2+1+\reductive^{1}_{\omega} 0.\]
But by Lemma \ref{LemmZero} we know that $\reductive^{1}_{\omega} 0=0$, so $\reductive^2_\omega\omega^{\varepsilon_0+1}=\omega\cdot 2+1$.

If instead we set $\vartheta=\varepsilon_0$, we obtain $N=N^2_\omega(\vartheta)=0$. Then, $\sigma_2[0]=0$ and \[\reductive^{\Theta[N]}_\Lambda(\varepsilon_0)= \reductive^{1}_\omega(e^{\omega} 1)=e 1=\omega,\]
i.e. $\reductive^2_\omega\varepsilon_0=\omega$.

We remark that, if $\uprho_\omega,\uprho$ represent the ranks with respect to $\mathcal I_\omega$ and $\mathcal I_1$, then
\[\uprho_\omega\omega^{\varepsilon_0+1}=
\le^\omega\omega^{\varepsilon_0+1}=0=\le (\omega\cdot 2+1)=
\uprho\reductive^2_\omega\omega^{\varepsilon_0+1},\]
while
\[\uprho_\omega\varepsilon_0=
\le^\omega \varepsilon_0=1=\le\omega =
\uprho\reductive^2_\omega\varepsilon_0.\]
\end{examp}

The above example suggests that $\reductive^2_\omega$ is rank-preserving. In Lemma \ref{partialisd} we will show that $\reductive^\Theta_\Lambda$ is always a $d$-map, so in general we always have $\uprho_\Lambda=\uprho\reductive^\Theta_\Lambda.$
In fact, these maps are $\Lambda$-reductive, but proving this will require several steps. Let us begin with a useful technical lemma.

\begin{lemm}\label{techdelt}
Given an ordinal $\Theta$ and a limit ordinal $\Lambda$,
\begin{enumerate}
\item $\reductive^\Theta_\Lambda[0,e^\Lambda\Theta]=[0,e\Theta]$
\item if $N>0$, $\reductive_\Lambda^\Theta \Delta^\Theta_\Lambda[N] = \Sigma_\Theta[N].$
\end{enumerate}
\end{lemm}

\proof
Assume both claims are true for $\Theta'$ when $\Theta'<\Theta$. Note that the first claim is trivial when $\Theta=0$ because then both sides of the equality are the singleton $\{0\}$, so we may assume $\Theta>0$. We shall begin by checking the inclusion
\[\reductive^\Theta_\Lambda \Delta^\Theta_\Lambda[N] \subseteq\Sigma_\Theta[N].\]
For $\vartheta\in \Delta^\Theta_\Lambda[N]$ we see that
\[\reductive^\Theta_\Lambda\vartheta= \sigma_\Theta[N]+ \reductive^{\Theta[N]}_{\Lambda}\le^{\lambda[N]}\vartheta \in\big[\sigma_\Theta[N],\sigma_\Theta[N+1]\big),\]
where the last step uses Claim 1 by induction on $\Theta[N]<\Theta$, given that $\le^{\Lambda[N]}\vartheta\in[0,e^{\Lambda}(\Theta[N])]$.

From this it easily follows that $\reductive^\Theta_\Lambda[0,e^\Lambda\Theta]\subseteq[0,e\Theta]$, since $\reductive^\Theta_\Lambda e^\Lambda\Theta=e\Theta\in[0,e\Theta]$, while for $\vartheta<e^\lambda\Theta$, we set $N=N^\Theta_\Lambda(\vartheta)$ and by Claim 2 see that
\[\reductive^\Theta_\Lambda\vartheta\in  \Sigma_\Theta[N]\subseteq[0,e\Theta].\]

Now let us check that
\[\Sigma_\Theta[N]\subseteq \reductive^\Theta_\Lambda \Delta^\Theta_\Lambda[N] .\]
Let $\xi=\sigma_\Theta[N] +\xi'$ with $\xi'\leq e({\Theta[N]})$. Then, using Claim 1 by induction on $\Theta[N]<\Theta$, $\xi'=\reductive^{\Theta[N]}_{\Lambda}\gamma'$ for some $\gamma'\leq e^{\Lambda}(\Theta[N])$; meanwhile, by Lemma \ref{deltalemm}, $\le^{\Lambda[N]}\Delta^\Theta_\Lambda[N]=[0,e^\Lambda(\Theta[N])]$, hence $\gamma'=\le^{\Lambda[N]}\gamma$ for some $\gamma \in\Delta^\Theta_\Lambda[N]$. It immediately follows that $\xi=\reductive^\Theta_\Lambda\gamma$.

To check the remaining inclusion of Claim 1, consider $\xi\in[0,e\Theta]$. If $\xi=e\Theta$, then evidently $\xi=\reductive^\Theta_\Lambda e^\Lambda(\Theta)$. Otherwise, $\xi\in[\sigma_\Theta[N], \sigma_\Theta[N+1])$ for some $N$. But by Claim 2,
\[\xi\in \reductive^\Theta_\Lambda\Delta^\Theta_\Lambda[N]\subseteq \reductive^\Theta_\Lambda[0,e^\Lambda\Theta],\]
as required.
\endproof

\begin{lemm}\label{partialisd}
Given countable ordinals $\Theta,\Lambda$ such that $\Lambda$ is infinite and additively indecomposable,
\[\reductive^\Theta_\Lambda:(e^{\Lambda}(\Theta)+1)_{\Lambda}\to (e(\Theta)+1)_1\]
is an onto $d$-map.
\end{lemm}

\proof
Assume the claim is true for all $\Theta'<\Theta$. Note that surjectivity is already proven in Lemma \ref{techdelt}. Also, the case $\Theta=0$ is trivial because then both sides of the equality are the singleton $\{0\}$, so that we may assume $\Theta>0$.

Pick $\vartheta \leq e^\Lambda(\Theta)$. If $\vartheta< e^\Lambda(\Theta)$, $\vartheta\in\Delta^\Theta_\Lambda[N]$ for some $N$, and since these sets are all open by Lemma \ref{deltalemm}, it suffices to observe that $\reductive^\Theta_\Lambda\upharpoonright\Delta^\Theta_\Lambda[N] $ is a $d$-map. But it is equal to $\reductive^{\Theta[N]}_{\Lambda}\le^{\Lambda[N]}$, which by induction on $\Theta[N]<\Theta$ and Lemma \ref{lisd} is a composition of $d$-maps. Hence it is open and continuous near $\vartheta$.

Otherwise, $\vartheta = e^\Lambda(\Theta)$. Here we claim that $\reductive^{\Theta}_\Lambda D^\Theta_\Lambda[N]=S_\Theta[N]$, from which openness and continuity are immediate.

We have that $D^\Theta_\Lambda[N]=\{e^\Lambda\Theta\}\cup\bigcup_{n>N}\Delta^\Theta_\Lambda[n],$ whereas $S_\Theta[N]=\{e\Theta\}\cup \bigcup_{n>N}\sigma_\Theta[N]$, and thus
\begin{align*}
\reductive^\Theta_\Lambda D^\Theta_\Lambda[N]&=\reductive^\Theta_\Lambda\{\vartheta\}\cup \bigcup_{n>N}\reductive^\Theta_\Lambda\Delta^\Theta_\Lambda[n]\\
&=\{e\Theta\}\cup \bigcup_{n>N}\Sigma_\Theta[N]\\
&=S_\Theta[N],
\end{align*}
where the second equality follows from Lemma \ref{techdelt}.

To check that $\reductive^\Theta_\Lambda$ is pointwise discrete, pick $\xi\leq e\Theta$. If $\xi=e\Theta$ then $(\reductive^\Theta_\Lambda)^{-1}\xi=\{e^\Lambda\Theta\}$, which is discrete.

Otherwise, $\xi\in\Sigma_\Theta[N]$ for some $N$ and thus $(\reductive^\Theta_\Lambda)^{-1}\xi=\Delta^\Theta_\Lambda[N] \cap(\reductive^\Theta_\Lambda)^{-1}\xi,$ which is discrete as $\reductive^\Theta_\Lambda\upharpoonright \Delta^\Theta_\Lambda[N]$ is a $d$-map.
\endproof

\begin{lemm}\label{thereismapa}
For any additively indecomposable limit ordinal $\Lambda$, the map $\reductive^\Theta_\Lambda$ is $\Lambda$-reductive.
\end{lemm}

\proof
By Lemma \ref{partialisd}, it suffices to prove that $\reductive^\Theta_\Lambda$ is $\lambda$-locally $\lambda$-determined for all $\lambda<1+\Lambda$; that is, that given $\vartheta<e^\Lambda \Theta $ there is a $\lambda$-neighborhood $U$ of $\vartheta$ with $\reductive^\Theta_\Lambda$ $\lambda$-determined on $U$.

Let $M$ be the largest natural such that $\Lambda[M]<\lambda$. We will consider two cases. First suppose that $\vartheta\in D^\Theta_\Lambda[M]$. Note that \[D^\Theta_\Lambda[M+1]=\big(e^{\Lambda[M]}\Theta[M] e^\Lambda\Theta\big]_{\Lambda[M]}\]
is $\lambda$-open. Let us check that $\reductive^\Theta_\Lambda$ is $\lambda$-determined on $D^\Theta_\Lambda[N+1]$. Assume that $\xi,\zeta\in D^\Theta_\Lambda[N+1]$ and $\le^\lambda\xi=\le^\lambda\zeta$. It readily follows that, for $n>M$, $\le^{\Lambda[n]}\xi=\le^{\Lambda[n]}\zeta$, and thus $n:=N^\Theta_\Lambda[\xi]=N^\Theta_\Lambda[\xi]$, from which we obtain
\[\reductive^\Theta_\Lambda\xi= \reductive^\Theta_\Lambda\zeta=\sigma_\Theta[n]+ \reductive^{\Theta[n]}_\Lambda \le^{\Lambda[n]}\xi.\] Thus we see that $\reductive^\Theta_\Lambda$ is $\lambda$-determined on the $\lambda$-open set $U=D^\Theta_\Lambda[M+1]$.

Otherwise we have that $\vartheta\in\Delta^\Theta_\Lambda[N]$ for some $N$ such that $\lambda\geq \Lambda[N]$. By induction on $\Theta[N]<\Theta$ there is a $(-\Lambda[N]+\lambda)$-neighborhood $U\subseteq [0,e^{\Lambda}(\Theta[N])]$ of $\le^{\Lambda[N]}\vartheta$ such that $\reductive^{\Theta[N]}_{\Lambda}$ is $(-\Lambda[N]+\lambda)$-determined on $U'$. By Lemma \ref{lisd}, $\le^{-\Lambda[N]}U'$ is $\lambda$-open, so that $U=\Delta^\Theta_\Lambda[N]\cap \le^{-\Lambda[N]}U'$ is $\lambda$-open as well. Moreover, $\reductive^{\Theta}_\Lambda$ is  ${\lambda}$-determined on $U$. To see this, assume that $\xi,\zeta\in U$ and $\le^{\lambda}\xi=\le^\lambda\zeta$; observe that the latter is equivalent to stating that \[\le^{-\Lambda[N]+\lambda}\le^{\Lambda[N]}\xi= \le^{-\Lambda[N]+\lambda} \le^{\Lambda[N]}\zeta.\]
It follows that
\begin{align}
\nonumber\reductive^\Theta_\Lambda\xi&= \sigma[N]+\reductive^{\Theta[N]}_\Lambda \le^{\Lambda[N]}\xi \\
\label{StepTwo} &= \sigma[N]+\reductive^{\Theta[N]}_\Lambda  \le^{\Lambda[N]}\zeta \\
 &=\nonumber
\reductive^\Theta_\Lambda\xi,
\end{align}
where \eqref{StepTwo} follows from the fact that $\le^{\Lambda[N]}\xi,\le^{\Lambda[N]}\zeta\in U'$ and $\reductive^{\Theta[N]}_\Lambda$ is $(-\Lambda[N]+\lambda)$-determined on $U'$.
\endproof

Let us conclude this section by extending $\reductive^\Theta_\Lambda$ to the case where $\Lambda$ is not necessarily additively indecomposable.

\begin{theorem}\label{thereismaptheo}
Given countable ordinals $\Theta,\Lambda$ with $\Lambda\geq 0$ there exists a $\Lambda$-reductive surjection
\[\reductive^\Theta_\Lambda:(e^{1+\Lambda}(\Theta)+1)_{1+\Lambda}\to(e(\Theta)+1)_1.\]
\end{theorem}

\proof
We may set $\reductive^\Theta_0={\sf id}$, and if $\Lambda+1$ is a successor, we may define $\reductive^\Theta_{\Lambda+1}=\le^{\Lambda}$ in view of Lemma \ref{LemmLogRed}.

If $\Lambda$ is infinite and additively indecomposable we have already defined $\reductive^\Theta_\Lambda$. Otherwise, if $\Lambda=\alpha+\omega^\beta$ with $\beta>1$ define $\reductive^\Theta_\Lambda= \reductive^\Theta_{\omega^\beta} \le^\alpha$, which is easily seen to be $\Lambda$-reductive.
\endproof


\section{Operations on ambiances}\label{secop}

In this section we shall review some operations, many of which were introduced in \cite{BeklGabel:2011}, that may be used to construct new provability ambiances from existing ones. First, let us observe that ambiances may be ``pulled back''.

If $\mathfrak X$ is an $(\alpha+\beta)$-BG space and $\mathfrak Y$ is a $\beta$-BG space, $f:\mathfrak X\to\mathfrak Y$ is an {\em $\alpha$-lift} if it is $\alpha$-reductive and, for all $\delta<\beta$, $f:\mathfrak X_{\alpha+\delta}\to\mathfrak Y_\delta$ is a $d$-map both with respect to the BG topologies and the shifted Icard topologies.

\begin{lemm}\label{pullalg}
Suppose that $\mathfrak X$ is a $(\xi+\zeta)$-BG-space and $\mathfrak Y$ is an idyllic $\zeta$-ambiance with algebra $\mathcal A$. Suppose further that $f:\mathfrak X\to \mathfrak Y$ is a $\xi$-lift. Then, $f^{-1}\mathcal A$ is an idyllic algebra on $\mathfrak X$.
\end{lemm}

\proof
To see that $d_{\xi+\delta}\upharpoonright f^{-1}\mathcal A=i_{1+\xi+\delta}\upharpoonright f^{-1}\mathcal A$, note that for $A\in\mathcal A$,
\[d_{\xi+\delta} f^{-1}A=f^{-1}d_{\xi+\delta} A=f^{-1}i_{1+\xi+\delta} A=i_{1+\xi+\delta} f^{-1}A.\]

That $d_\delta\upharpoonright f^{-1}\mathcal A=i_{1+\delta}\upharpoonright f^{-1}\mathcal A$ for $\delta<\alpha$ follows from Lemma \ref{absolutely} since $f^{-1}A$ is always $\delta$-absolute.
\endproof

One of the basic ways of including topological spaces into a larger space is by their {\em topological sum}. In the case of ordinal spaces, this is closely tied to the sum of ordinals, as has already been observed in \cite{BeklGabel:2011}:

\begin{definition}
Given families of sets $\mathcal A\subseteq\mathcal P(\Xi)$ and $ \mathcal B\subseteq\mathcal P(\Theta)$, where $\Xi,\Theta$ are ordinals, we define $\mathcal A\oplus\mathcal B$ to be the family of subsets of $\Xi+\Theta$ of the form
\[S=S_0\cup (\eta+S_1),\]
with $S_0\in\mathcal A$ and $S_1\in\mathcal B$.
\end{definition}

The following lemma is standard and easy to check:

\begin{lemm}
If $\Xi,\Theta$ are ordinals and $\mathcal A\subseteq\mathcal P(\Xi),\mathcal B\subseteq\mathcal P(\Theta)$ are topologies, then $\mathcal A\oplus\mathcal B$ is a topology on $\Xi+\Theta$.

If $\Xi$ is a successor and both topologies are Icard or BG, then $\mathcal A\oplus\mathcal B$ is Icard or BG, respectively.
\end{lemm}

In view of this we define, given $\Lambda$-polytopologies $\mathfrak X=\langle \Xi,\vec{\mathcal T}\rangle$ and $\mathfrak Y=\langle \Theta,\vec{\mathcal S}\rangle$, the sum
\[\mathfrak X\oplus\mathfrak Y=\langle \Xi+\Theta,\langle\mathcal T_\lambda\oplus\mathcal S_\lambda\rangle_{\lambda<\Lambda}\rangle. \]

We may also apply the sum operation to $d$-algebras:

\begin{lemm}
Supppose $\mathfrak X=\langle\Xi+1,\vec{\mathcal T},\mathcal A\rangle$ and $\mathfrak Y=\langle \Theta,\vec{\mathcal S},\mathcal A\rangle$ are ambiances. Then, $\mathfrak X\oplus\mathfrak Y$ equipped with $\mathcal A\oplus\mathcal B$ is also an ambiance. Further, if both ambiances are idyllic, then so is the corresponding sum.
\end{lemm}

We shall not present a proof, as this result is fairly obvious once we observe that the all relevant operations may be carried out independently within the two disconnected subspaces.

The last major topological construction needed for the completeness proof is the notion of a {\em $d$-product}. Let $\mathfrak X=\langle [0,\Xi],\langle\mathcal T_\lambda\rangle_{\lambda<\Lambda}\rangle$ and $\mathfrak Y=\langle [0,\Theta],\langle\mathcal S_\lambda\rangle_{\lambda<\Lambda}\rangle$ be polytopologies and define $\Xi\otimes_d\Theta=-1+(1+\Xi)(1+\Theta)$.

Let $G_1$ be the set of those $\xi\leq\Xi\otimes_d\Theta$ of the form $(1+\Xi)\alpha$ with $\alpha\in\sf Lim$ and $G_0$ be its complement; we will call $G_0,G_1$ the {\em components} of $[0,\Xi\otimes_d\Theta]$. Every $\xi\in G_0$ can be written uniquely in the form \[\xi=-1+(1+\Xi)\alpha+(1+\xi_0)\]
with $\xi_0\leq \Xi$; to see this, write $\xi=-1+(1+\Xi)\gamma+\delta$ with $\delta<1+\Xi$. If $\delta>0$ we may take $\xi_0=-1+\delta$. If $\delta=0$, write $\gamma=\omega\cdot\gamma' +n$. If $n=0$ we would have $\xi\in G_1$, so $n>0$ and $\xi=-1+(1+\Xi)(\omega\cdot \gamma'+(n-1))+(1+\Xi)$. We may then take $\xi_0=\Xi$. Note that the term $-1$ only affects the expression when both $\alpha$ and $\Xi$ are finite, in which case it becomes $\xi=(1+\Xi)\alpha+\xi_0$.

With this we define $\pi_0:G_0\to [0,\Xi]$ by setting $\pi_0(\xi)=\xi_0$. Similarly, define $\pi_1:[0,\Xi\otimes_d\Theta]\to [0,\Theta]$ as follows:
\begin{itemize}
\item if $\xi=-1+(1+\Xi)\alpha+1+\pi_0\xi\in G_0$ then $\pi_1\xi=-1+\alpha+1$,
\item if $\xi=(1+\Xi)\alpha\in G_1$ then $\pi_1\xi=\alpha$.
\end{itemize}

\begin{examp}
Let $\Xi=\omega^2$ and $\Theta=\omega+1$. Then, $\Xi\otimes_d\Theta=\omega^3+\omega^2$. Elements of $G_1$ are those of the form $(1+\Xi)\alpha$ with $\alpha\in \sf Lim$, i.e. those of the form $\omega^3\gamma$.

Let us compute some projections of specific elements.
\begin{itemize}
\item $3\not\in G_1$, and $3=-1+(1+\omega^2)0+1+3$, so $\pi_0 3=3$ and $\pi_1 3=-1+0+1=0$.
\item $\omega^2+3=-1+(1+\omega^2)1+1+2$ and hence $\pi_0(\omega^2+\omega)=2$ whereas $\pi_1\omega=-1+1+1=1$.
\item $\omega^3+\omega+3=-1+(1+\omega^2)\omega+1+\omega+3$, so we have $\pi_0 (\omega^3+\omega+3)=\omega+3$ and $\pi_1(\omega^3+\omega+3)=-1+\omega+1=\omega+1$.
\item $\omega^3=(1+\omega^2)\omega$ is the only element of $G_1$, so $\pi_0(\omega^3)$ is undefined but $\pi_0(\omega^3)=\omega$.
\end{itemize}
\end{examp}

Although the projections are not injective, they have natural bijective restrictions.

\begin{lemm}\label{LemmBiject}
Let $\Theta,\Xi$ be ordinals. The restrictions $\pi_0:\pi^{-1}_1(\alpha)\to \Xi+1$ with $\alpha\in [0,\Theta]\setminus \sf Lim$ and $\pi_1:G_1\to [0,\Theta]\cap \sf Lim$ are strictly increasing and onto.
\end{lemm}

\proof
Pick $\alpha\in[0,\Theta]\setminus\sf Lim$. Let $\beta=\alpha$ if $\alpha$ is finite or $\beta+1=\alpha$ if $\alpha$ is infinite, and let $\gamma=-1+\beta+1$. By observation on the definition of $\pi_1$, $\pi^{-1}_1\{\alpha\}$ is an interval of the form $J_\alpha=[\gamma,\gamma+\Xi]$. On this interval, $\pi_0(\xi)=-\gamma+\xi$, which is clearly increasing and onto $[0,\Xi]$.

For the second claim, given $\xi\in[0,\Xi\otimes_d\Theta]$, $\xi\in G_1$ if and only if $\pi_1\xi\in \sf Lim$, and indeed $\pi_1\upharpoonright G_1$ is given by the map $(1+\Xi)\alpha\mapsto \alpha$ which is increasing and onto $[0,\Theta]\cap\sf Lim$.
\endproof

It will also be useful to see how the projections treat hyperlogarithms.

\begin{lemm}\label{LemmLogProy}
Let $\Xi,\Theta$ be ordinals, and write $1+\Xi=\omega^\alpha+\beta$ with $\beta<1+\Xi$. Let $G_0,G_1$ be the components of $[0,\Xi\otimes_d\Theta]$.

Then, for $\xi\in G_0$ we have that $\le \pi_0\xi=\le\xi$, whereas for $\xi\in G_1$ we have that $\le\xi=\alpha+\le\pi_1\xi$.

Moreover, for $\lambda>1$, we have that $\le^\lambda \pi_0\xi=\le^\lambda\xi$, whereas for $\xi\in G_1$ we have that $\le^\lambda\xi=\le^\lambda\pi_1\xi$.
\end{lemm}

\proof
Note that the claim for $\lambda>1$ is an immediate consequence of the first claim, so we focus on $\lambda=1$.

An ordinal $\xi\in G_0$ is of the form $-1+(1+\Xi)\alpha+1+\pi_0\xi$, and hence $\le\xi=\le\pi_0\xi$; note that this holds even when $\pi_0\xi=0$, in which case $\le\xi=\le\pi_0\xi=0$. Any ordinal $\xi\in G_1$ is of the form $(1+\Xi)\omega(1+\zeta)$, with $\omega(1+\zeta)\leq 1+\Theta$. Write $1+\zeta=\gamma+\omega^\delta$. Then, we see that
\begin{equation}\label{decomposition}
\xi=(\omega^\alpha+\beta)\omega(\gamma+\omega^\delta)= \omega^{\alpha+1}\gamma+\omega^{\alpha+1+\delta}.
\end{equation}
Hence $\le\xi=\alpha+1+\delta$. However, $\pi_1\xi=\omega\gamma+\omega^{1+\delta}$, and thus $\le\pi_1\xi=1+\delta$. The result follows.
\endproof

We are now ready to define the $d$-product of spaces:

\begin{definition}[$d$-product]\label{dprod}
Let $\mathfrak X=\langle \Xi+1,\langle\mathcal T_\lambda\rangle_{\lambda<\Lambda}\rangle$ and $\mathfrak Y=\langle \Theta+1,\langle\mathcal S_\lambda\rangle_{\lambda<\Lambda}\rangle$ be polytopologies.

 For $\lambda<\Lambda$, define a topology $\mathcal O_\lambda$ on $[0,\Xi\otimes_d\Theta]$ to be generated by sets of the forms
\begin{itemize}
\item $\pi^{-1}_0 U \cap \pi^{-1}_1 \{\alpha\}$, where $U\subseteq [0,\Xi]$ is $\lambda$-open and $\alpha\leq \eta$ is finite or a successor, or
\item $\pi_1^{-1} U$, where $U\subseteq [0,\Theta]$ is $\lambda$-open.
\end{itemize}

We denote the resulting space $\langle[0,\Xi\otimes_d\Theta],\langle\mathcal O_\lambda\rangle_{\lambda<\Lambda}\rangle$ by $\mathfrak X\otimes_d\mathfrak Y$.
\end{definition}

This definition will be sufficient for our purposes, but the $d$-product of spaces is treated with much more generality and detail in \cite{BeklGabel:2011}, which gives it a slightly different presentation. There, the following properties are established:

\begin{lemm}\label{LemmDprodTop}
If $\mathfrak X=\langle[0,\Xi],\langle\mathcal T_\lambda\rangle_{\lambda<\Lambda}\rangle$ and $\mathfrak Y=\langle[0,\Theta],\langle\mathcal S_\lambda\rangle_{\lambda<\Lambda}\rangle$ are regular polytopologies and $\mathfrak Z=\mathfrak X\otimes_d\mathfrak Y$, then
\begin{enumerate}
\item $G_0$ is $0$-open and $G_1$ is $1$-open,
\item for every $\xi\in [0,\Xi]$, $\pi_0^{-1}\xi$ is $0$-dense in $G_1$,\label{LemmDprodTop2}
\item $A\subseteq G_1$ is open in $ (\mathcal T_0\otimes_d\mathcal S_0)\upharpoonright G_1$ (i.e., in the subspace topology) if and only if $G_0\cup A$ is open in $\mathcal T_0\otimes_d\mathcal S_0$.\label{LemmDprodTop3}
\end{enumerate}
\end{lemm}

\proof
It is easy to see that $G_0$ is $0$-open, since \[G_0=\bigcup_{\alpha\in [0,\vartheta]\setminus{\sf Lim}}\pi^{-1}_0[0,\eta]\cap \pi^{-1}_1\{\alpha\}.\]
$G_1$ is $1$-open since $G_1=\pi_1^{-1}(1,\Theta]_1$.

To see that $\pi_0^{-1}\xi$ is $0$-dense in $G_1$, pick $\zeta\in G_1$. A basic neighborhood of $\zeta$ is of the form $\pi^{-1}_1 U$ with $U\subseteq [0,\Theta]$ $0$-open. By Lemma \ref{rankopen}, $U$ contains an element $\delta$ with rank $0$, i.e. zero or a successor. Then, by Lemma \ref{LemmBiject}, $\pi_0\pi^{-1}_1\delta=[0,\Xi]$, i.e. $\xi\in \pi_0\pi^{-1}_1\delta$ or, equivalently, $\varnothing\not=\pi^{-1}_0\xi\cap \delta\subseteq \pi^{-1}_1 U$. Since $B$ was arbitrary, the result follows.

For the third claim, note that $A$ is $0$-open in $G_1$ if and only if there is a $0$-open set $U$ such that $U\cap G_1=A$. But then we have that $G_0\cup U= G_0\cup (G_0\cap U)\cup (G_1\cap U)=G_0\cup A$, and since it is a union of opens it must be open. Conversely, if $G_0\cup A$ is $0$-open, then $(G_0\cup A)\cap G_1=A$ and hence $A$ is $0$-open in the subspace topology.

\begin{lemm}\label{LemDprodMap}
If $\mathfrak X=\langle[0,\eta],\langle\mathcal T_\lambda\rangle_{\lambda<\Lambda}\rangle$ and $\mathfrak Y=\langle[0,\vartheta],\langle\mathcal S_\lambda\rangle_{\lambda<\Lambda}\rangle$ are regular polytopologies and $\mathfrak Z=\mathfrak X\otimes_d\mathfrak Y$, then
\begin{enumerate}
\item for all $\lambda<\Lambda$,
$\pi_0\colon (\mathfrak Z_\lambda\upharpoonright G_0)\to \mathfrak X_\lambda$ is a $d$-map.
\item for all $\lambda<\Lambda$ and $\alpha\in [0,\Theta]\setminus\sf Lim$,
$\pi_0\colon (\mathfrak Z_\lambda\upharpoonright \pi_1^{-1}\alpha)\to \mathfrak X_\lambda$ is a homeomorphism.
\item $\pi_1$ is $\lambda$-continuous and $\lambda$-open,
\item $\pi_1\upharpoonright G_1\colon (\mathfrak Z_\lambda\upharpoonright G_1)\rightarrow (\mathfrak Y_\lambda\upharpoonright {\sf Lim})$ is a homeomorphism.
\end{enumerate}
\end{lemm}

\proof
Items 1 and 3 are proven by showing that images of basic opens are open, as are preimages of basic opens. We will not provide the details, but only show as an example that $\pi_0\colon (\mathfrak Z_\lambda\upharpoonright G_0)\to \mathfrak X_\lambda$ is open. The $\lambda$-topology on $\mathfrak Z$ is generated by sets of the form $\pi_0^{-1}V\cap \pi^{-1}_1\alpha$ with $\alpha\in[0,\Theta]\setminus \sf Lim$ and $V\subseteq [0,\Xi]$ $\lambda$-open, or $\pi^{-1}_1 U$ with $U\subseteq [0,\Theta]$ $\lambda$-open. By Lemma \ref{LemmBiject}, $\pi_0\upharpoonright \pi^{-1}_1\alpha$ is a bijection and hence $\pi_0(\pi_0^{-1}V\cap \pi^{-1}_1\alpha)=V$, which is open. In the second case, $U$ either contains some $\alpha\in [0,\Theta]\setminus \sf Lim$, or it does not. If it does, $\pi_0 (U\cap G_0)=[0,\Theta]$, and if not, $\pi_0 (U\cap G_0)=\varnothing$, both of which are open. Note that the fact that $\pi_0$ is pointwise discrete follows also by the injectivity of $\pi_0\upharpoonright \pi^{-1}_1\alpha$, as given $\zeta\in\pi^{-1}_0\xi$, $\zeta$ is the only point of $\pi^{-1}_0\xi$ contained in the open set $\pi^{-1}_0[0,\Xi]\cap \pi^{-1}_1\pi_1\zeta$.

Items 2 and 4 follow from Items 2 and 3, respectively, and Lemma \ref{LemmBiject}; $\pi_0\upharpoonright \pi^{-1}_1\alpha$ is a restriction of a continuous and open map to an open set, hence it is continuous and open, and as it is bijective, it is a homeomorphism. Since $G_1=\pi^{-1}_1([0,\Theta]\cap{\sf Lim})$, we have that $\pi_1\upharpoonright G_1=\pi_1\upharpoonright\pi^{-1}_1 ([0,\Theta]\cap{\sf Lim})$ is continuous and open\footnote{In general, if $f\colon X\to Y$ is a continuous (open) function and $A\subseteq Y$, then $f\upharpoonright f^{-1}A$ is continuous (open).} and hence, being bijective, a homeomorphism.
\endproof

\begin{lemm}
Suppose that $\mathfrak X$ and $\mathfrak Y$ are regular polytopologies and $\mathfrak Z=\mathfrak X\otimes_d\mathfrak Y$. If $\mathfrak X,\mathfrak Y$ are both BG spaces, then so is $\mathfrak Z$. Similarly, if $\mathfrak X,\mathfrak Y$ are shifted Icard spaces, then so is $\mathfrak Z$.
\end{lemm}

\proof
\cite[Lemma 7.6]{BeklGabel:2011} states that, if both $\mathfrak X'$ and $\mathfrak Y'$ are limit-maximal extensions of $\mathfrak X$ and $\mathfrak Y$, respectively, then $\mathfrak X'\otimes_d\mathfrak Y'$ is a limit-maximal extension of $\mathfrak X\otimes_d\mathfrak Y$. Hence in order to show that the $d$-product of BG-spaces is BG, it suffices to show that the $d$-product of their underlying Icard space is Icard.

So let $\Xi,\Theta,\Lambda$ be ordinals, $\Omega=\Xi\otimes_d\Theta$ and $G_0,G_1$ be the components of $\Omega+1$ with projections $\pi_0,\pi_1$. Let $\mathcal O_\lambda=\mathcal I^{\Xi+1}_{1+\lambda}\otimes_d\mathcal I^{\Theta+1}_{1+\lambda}$ for each $\lambda<\Lambda$. We claim that $\mathcal O_\lambda=\mathcal I^{\Omega+1}_{1+\lambda}$ for all $\lambda<\Lambda$.

We proceed by induction on $\lambda$. The base case, when $\lambda=1$, is proven in \cite[Section 7]{BeklGabel:2011}. For $\lambda=\eta+1$ a successor, we have that $\mathcal I_{1+\lambda}=\dot{\mathcal I}_{1+\eta}$. Thus we must show that $\dot{\mathcal O}_{\eta}=\dot{\mathcal I}^{\Xi+1}_{1+\eta}\otimes_d\dot{\mathcal I}^{\Theta+1}_{1+\eta}$, under the hypothesis that $\mathcal O_\eta=\mathcal I_{1+\eta}^{ \Omega+1}$. To show that $\dot{\mathcal O}_{\eta}\subseteq\dot{\mathcal T}_\eta\otimes_d\dot{\mathcal S_\eta}$, it suffices to show that $(\delta,\Omega]_{1+\eta}\in \dot{\mathcal I}^{\Xi+1}_{1+\eta}\otimes_d\dot{\mathcal I}^{\Theta+1}_{1+\eta}$ for all $\delta$.

First assume that $\xi\in G_0\cap (\delta,\Omega]_{1+\eta}$. Then, by Lemma \ref{LemmLogProy} we have that $\le^{1+\eta}\xi=\le^{1+\eta}\pi_0\xi$ and thus
\[\xi\in \pi_0^{-1}(\delta,\Xi]_{1+\eta}\cap \pi^{-1}_0\pi_1\xi\subseteq (\delta,\Omega+1]_{1+\eta},\]
and $\pi_0^{-1}(\delta,\Xi]\cap\pi^{-1}_0\pi_1\xi\in \dot{\mathcal T}_\eta\otimes_d\dot{\mathcal S_\eta}$. If $\xi \in G_1$, we consider two subcases. If $\eta=0$ then write $1+\Xi=\omega^\alpha+\beta$. Without loss of generality we may assume that $\delta\geq \alpha$, for otherwise we have that $(\alpha,\Omega+1]_{1}\subseteq(\delta,\Omega+1]_{1}$. Then, we have by Lemma \ref{LemmLogProy} that $(\delta,\Omega+1]_{1}=\pi_1^{-1}(-\alpha+\delta,\Theta+1]_1$. The case where $\eta>0$ is similar, but here we have simply that \[\xi\in G_1\cap \pi_1^{-1}(\delta,\Theta+1]_{1+\eta}\subseteq(\delta,\Omega+1]_{1+\eta}.\] Since $\Omega+1=G_0\cup G_1$, we conclude that $(\delta,\Omega+1]_{1+\eta}\in \dot{\mathcal I}^{\Xi+1}_{1+\eta}\otimes_d\dot{\mathcal I}^{\Theta+1}_{1+\eta}$.

For the other inclusion, we note that $\dot{\mathcal I}^{\Xi+1}_{1+\eta}\otimes_d\dot{\mathcal I}^{\Theta+1}_{1+\eta}$ is generated by ${\mathcal I}^{\Xi+1}_{1+\eta}\otimes_d{\mathcal I}^{\Theta+1}_{1+\eta}$ and sets of the form $\pi^{-1}_0(\delta,\Xi]_{1+\eta}\cap \pi^{-1}_1\vartheta$ for $\vartheta\in[0,\Theta]\setminus\sf Lim$, or of the form $\pi^{-1}_1(\delta,\Theta]_{1+\eta}$. It remains to check that they are both open in $\dot{\mathcal O}_\eta$. But in the first case we have that
\[\pi^{-1}_0(\delta,\Xi]_{1+\eta}\cap \pi^{-1}_1\vartheta=(\delta,\Omega]_{1+\eta}\cap \pi^{-1}_1\vartheta, \]
which is an intersection of opens in $\dot{\mathcal O}_{\eta}$ and hence open, whereas for $\xi\in\pi^{-1}_1(\delta,\Theta]_{1+\eta}$ we have that, if $\xi\in G_0$, then $\pi^{-1}_0\pi_1\xi \subseteq \pi^{-1}_1(\delta,\Theta]_{1+\eta}$ and $\pi^{-1}_0\pi_1\xi\in \mathcal O_\eta$, whereas if $\xi\in G_1$, then
\[G_1\cap \pi^{-1}_1(\gamma+\delta,\Omega]_{1+\eta} \subseteq \pi^{-1}_1(\delta,\Theta]_{1+\eta},\]
where $\gamma=\alpha$ if $\eta=0$ and $\gamma=0$ if $\eta>0$.

Finally, we consider the case where $\lambda\in \sf Lim$. But here it is very easy to check that
\[\mathcal O_\lambda=\bigsqcup_{\eta<\lambda}{\mathcal I}^{\Xi+1}_{\eta}\otimes_d{\mathcal I}^{\Theta+1}_{\eta}\stackrel{\rm IH}=\bigsqcup_{\eta<\lambda}\mathcal I_\eta^{ \Omega+1}=\mathcal I_\lambda^{\Omega+1}. \]
\endproof

With these ingredients, we define the {\em $d$-product} of algebras: 

\begin{definition}[$d$-product of algebras]\label{dprodal}
Given $d$-algebras $\mathcal A,\mathcal B$ based on $\Xi+1$, $\Theta+1$ and letting $G_0,G_1$ be the components of $[0,\Xi\otimes_d\Theta]$, we define $\mathcal A\otimes_d\mathcal B$ to be the algebra of all sets $S$ of the form
\[S=\pi_0^{-1}(S_0)\cup \pi^{-1}_1(S_1\cap G_1),\]
where $S_0\in \mathcal A$ and $S_1\in \mathcal B$.
\end{definition}

Of course, we would like for the $d$-product of algebras to be itself a $d$-algebra. The next lemma will be useful in showing this.

\begin{lemm}\label{dmult}
Let $\mathfrak X,\mathfrak Y$ be polytopologies based on ordinals $\Xi+1,\Theta+1$, respectively. Let $d_\xi$ denote the $\xi$-derived set operator on $\mathfrak Y$ and $d'_\xi$ on $\mathfrak X\otimes_d\mathfrak Y$.

Then, for any $E\subseteq i_1[0,\Theta]$ and $\lambda<\Lambda$, $d'_\lambda \pi^{-1}_1E=\pi^{-1}_1d_\lambda E$.
\end{lemm}

\proof
Assume first that $\xi\in d'_\lambda \pi^{-1}_1( E)$. Note that this immediately implies that $\xi\not\in G_0$, since the latter is $0$-open. 

Then, for every $\lambda$-neighborhood $U$ of $\xi$ there is $\zeta\not=\xi\in \pi^{-1}_1( E)\cap U$. Now, if $V$ is a $\lambda$-neighborhood of $\pi_1 \xi$ in $\Theta+1$, then $\pi^{-1}_1( V)$ is a $\lambda$-neighborhood of $\xi$, so that there is $\zeta\not=\xi\in \pi^{-1}_1(E)\cap  \pi^{-1}_1 (V)$, hence $\pi_1\zeta\in E\cap V$. Further, by Lemma \ref{LemmBiject}, $\pi_1$ is injective on $G_1=\pi_1^{-1}(0,\Theta]_1$ so $\pi_1\zeta\not=\pi_1\xi$. Since $V$ was arbitrary, it follows that $\pi_1\xi\in d_\lambda E$.

Conversely, if $\pi_1\xi \in d_\lambda E$, then every $\lambda$-neighborhood $V$ of $\pi_1\xi$ contains some $\zeta\not=\xi\in E$. Consider a $\lambda$-neighborhood $U=\pi^{-1}_1V$ of $\xi$, where $V\subseteq [0,\Theta]$ is $\lambda$-open. Since $\pi_1\xi \in d_\lambda E$ there is some $\zeta\in E\cap V$ with $\zeta\not=\xi$ and since $\pi_1$ is onto we have that $\zeta=\pi_1\zeta'$ for some (unique) $\zeta'\in G_1$ and thus $\zeta'\in U\cap \pi^{-1}_1 E$. Since $U$ was arbitrary and $\zeta'\not=\xi$, we conclude that $\xi\in d'_\lambda\pi^{-1}_1E$.
\endproof

\begin{lemm}\label{prodil}
Supppose $\mathfrak X=\langle[0,\Xi],\vec{\mathcal T},\mathcal A\rangle$ and $\mathfrak Y=\langle[0,\Theta],\vec{\mathcal S},\mathcal B\rangle$ are ambiances. Then, $\mathfrak X\otimes_d\mathfrak Y$ equipped with $\mathcal A\otimes_d\mathcal B$ is also an ambiance.  Further, if both ambiances are idyllic, then so is the corresponding product.
\end{lemm}

\proof
Let $S=\pi_0^{-1}(S_0)\cup \pi^{-1}_1(S_1\cap G_1)\in\mathcal A\otimes_d\mathcal B$ and $\lambda<\Lambda$. Let us check that $d_\lambda S\in\mathcal A\otimes_d\mathcal B$.

If $\lambda=0$ and $S_0\not=\varnothing$, then by Lemma \ref{LemmDprodTop}.\ref{LemmDprodTop2}, $\pi^{-1}_0S_0$ is $0$-dense in $G_1$ and thus $G_1\subseteq d_0 S$.

Meanwhile, $G_0$ is $0$-open and $\pi_0$ is a $d$-map so for $\xi\in G_0$, $\xi\in d_\lambda S$ if and only if $\pi_0\xi\in d_\lambda S_0$, and we conclude that
\[d_0 S=\pi_0^{-1}d_0S_0\cup G_1\in \mathcal A\otimes_d\mathcal B.\] 
By similar reasoning,
\[i_1 S=\pi_0^{-1}i_1S_0\cup G_1\in \mathcal A\otimes_d\mathcal B,\] 
and if the original structures are idyllic this is evidently equal to $d_0 S$ as well.

Now suppose $\lambda=0$ and $S_0=\varnothing$. In this case, $S=\pi^{-1}_1 (S_1\cap G_1)$ and, by Lemma \ref{dmult}, $d_\lambda S=\pi^{-1}_1 d_\lambda S_1\in\mathcal A\otimes_d\mathcal B$. If the original algebras were idyllic, we note further that
\[d_0 S=\pi^{-1}_1 d_0 S_1=\pi^{-1}_1 i_{1} S_1=i_{1}S.\]

Finally, if $\lambda>0$, then both projections are $d$-maps with respect to the $\lambda$-topology and
\[d_\lambda S=d_\lambda \pi_0^{-1} S_0\cup d_\lambda\pi^{-1}_1 (S_1\cap G_1)=\pi_0^{-1}d_\lambda  S_0\cup \pi^{-1}_1d_\lambda (S_1\cap G_1),\]
with the analogous equalities holding for $i_{1+\lambda}S$, from which all required claims follow easily.
\endproof

With this we conclude the topological constructions we shall need. Now, we turn to the last ingredient in the completeness proof: the modal logic $\sf J$.


\section{The logic $\sf J$}\label{jsec}
As we have seen, $\glp_\Lambda$ has no non-trivial Kripke frames. In order to work around this issue, we pass to a weaker logic, Beklemishev's ${\sf J}$. This was introduced in \cite{levj} and here we only review the necessary results without proof. For this logic we shall only use modalities $n<\omega$ and replace Axiom \ref{glpfour} of $\glp_\Lambda$ by the two axioms
\begin{enumerate}
\item[6.] $[n]\phi\to[m][n]\phi$, for $n\leq m$ and
\item[7.] $[n]\phi\to[n][m]\phi$, for $n< m$.
\end{enumerate}

The logic ${\sf J}$ is sound and complete for the class of finite Kripke models $\langle W,\langle<_n\rangle_{n<N},\lb\cdot\rb\rangle$ such that
\begin{enumerate}
\item the relations $<_n$ are transitive and well-founded,
\item if $n<m$ and $w<_m v$ then $\mathop <_n(w)=\mathop <_n(v)$ and
\item if $n<m$ then $w<_m v<_n u$ implies that $w<_n u$.
\end{enumerate}
Here, $\mathop <_n(w)=\{v:v<_n w\}$. It will also be convenient to define $w\ll_n v$ if for some $m\geq n$, $w<_m v$. Let $\sim_n$ denote the symmetric, transitive, reflexive closure of $\ll_n$ and let $[w]_n$ denote the equivalence class of $w$ under $\sim_n$. Define $[w]_{n+1}<_n[v]_{n+1}$ if there exist $w'\in[w]_{n+1}, v'\in [v]_{n+1}$ such that $w'<_n v'$.

Then, say $W$ is {\em tree-like} if
\begin{enumerate}
\item for each $w\in W$ and $n\leq N$, $[w]_n/\sim_{n+1}$ is a tree under $<_n$ and
\item if $[w]_{n+1}<_n[v]_{n+1}$ then $w<_n v$.
\end{enumerate}


With this we may state the following completeness result from \cite{levj}:

\begin{lemm}\label{jcomp}
Any $\sf J$-consistent formula can be satisfied on a finite, tree-like $\sf J$-frame.
\end{lemm}

Thus if we can reduce $\glp_\Lambda$ to $\sf J$, we will immediately obtain finite Kripke models. For this, given a formula $\phi$, let $N$ be the largest modality appearing in $\phi$ and define
\[M(\phi)=\bigwedge_{\substack{[n]\psi\in{\rm sub}(\phi)\\
n< m\leq N}}
[n]\psi\to [m]\psi.
\]
Then we set $M^+(\phi)=M(\phi)\wedge\bigwedge_{n\leq N}[n]M(\phi)$.

The following is also proven in \cite{levj}:
\begin{lemm}\label{glptoj}
For any formula $\phi\in {\sf L}_\omega$, $\glp_\omega\vdash\phi$ if and only if
\[{\sf J}\vdash M^+(\phi)\to \phi.\]
\end{lemm}

To prove completeness, it then suffices to construct a $\sf J$-model of a given formula and then ``pull back'' the valuations onto a topological model. Hence it is important to identify the appropriate maps for such pullbacks.

First observe that a partially ordered set $\langle W,<\rangle$ can be identified with a topological space by letting $U\subseteq W$ be open if, whenever $v<w$ and $w\in U$, it follows that $v\in U$. If $\mathfrak W=\langle W,\langle <_n\rangle_{n\leq N}\rangle$ is a $\sf J$-frame, we will let $\mathfrak W_n$ be the topological space associated to $\langle W,<_n\rangle$. Below, we say $w\in W$ is a {\em hereditary $n$-root} if $w$ is $<_k$-maximal for all $k\geq n$.

\begin{definition}\label{jmap}
Let $\mathfrak X=\langle X,\langle\mathcal T_n\rangle_{n\leq N}\rangle$ be a polytopological space and $\mathfrak W=\langle W,\langle <_n\rangle_{n\leq N}\rangle$ a $\sf J$-frame.

A function $f:\mathfrak X\to\mathfrak W$ is a {\em $\sf J$-map} if
\begin{enumerate}
\item $f:\mathfrak X_N\to \mathfrak W_N$ is a $d$-map
\item $f:\mathfrak X_n\to\mathfrak  W_n$ is open for $n<N$
\item if $n<N$ and $w$ is a hereditary $(n+1)$-root then $f^{-1}(\mathop \ll_n(w))$ and \[f^{-1}(\mathop {\ll_n}(w)\cup \{w\})\]
are $n$-open and
\item if $n<N$ and $w$ is a hereditary $(n+1)$-root then $f^{-1}(w)$ is $n$-discrete.
\end{enumerate}
\end{definition}

\begin{examp}
Consider a simple $\sf J$-frame $\mathfrak W$ with three worlds, $u,v,w$, such that $u<_0 w$ and $v<_0 w$. Let us define a $\sf J$-map $f\colon (\omega+1)_1\rightarrow \mathfrak W$.

The worlds $u,v$ are isolated, and hence we must have that $\xi$ is isolated whenever $f(\xi)\in\{u,v\}$. It follows that $f^{-1}\{u,v\}=[0,\omega)$. Meanwhile, the only possible value for $f(\omega)$ is $w$ as it is the only point that is not isolated and $d$-maps preserve rank.

Now, notice that any neigbhorhood of $w$ contains both $u$ and $v$. Thus we will need for every neighborhood of $\omega$ to intersect both $f^{-1}(u)$ and $f^{-1}(v)$. A simple way to achieve this is to let $f^{-1}(u)$ be the set of even numbers and $f^{-1}(v)$ the set of odds. Thus we may define
\[f(\xi)=
\begin{cases}
u&\text{if $\xi=2n<\omega$,}\\
v&\text{if $\xi=2n+1<\omega$,}\\
w&\text{if $\xi=\omega$.}
\end{cases}
\]
One can then easily check that the function $f$ thus defined is in fact a $\sf J$-map.
\end{examp}

With this we have the following, proven in \cite{BeklGabel:2011}:

\begin{lemm}\label{truthpres}
If $\mathfrak W=\langle W,\langle<_n\rangle_{n\leq N},\lb\cdot\rb\rangle$ is a $\sf J$-model such that $M^+(\phi)$ is valid on $\mathfrak W$, $\mathfrak X=\langle X,\langle\mathcal T_n\rangle_{n\leq N}\rangle$ is a GLP-space and $f:\mathfrak X\to\mathfrak W$ a ${\sf J}$-map, then there is a valuation $\lp\cdot\rp$ on $\mathfrak X$ such that $\lp\psi\rp=f^{-1}\lb\psi\rb$ for all $\psi\in{\rm sub}(\phi)$.
\end{lemm}

We conclude with a simple observation, also established in \cite{BeklGabel:2011}:

\begin{lemm}\label{dj}
Let $\mathfrak X,\mathfrak Y$ be polytopological spaces and $\mathfrak W$ a $\sf J$-frame.

Then, if $f:\mathfrak X\to\mathfrak Y$ is a $d$-map and $g:\mathfrak Y\to\mathfrak W$ is a $\sf J$-map it follows that $gf$ is a $\sf J$-map.
\end{lemm}

In the next section we shall exploit the completeness of the logic $\sf J$ for finite frames together with Lemma \ref{truthpres} to construct GLP-ambiances satisfying any consistent formula.


\section{Completeness}\label{seccomp}

Given an increasing sequence of ordinals $\vec\lambda=\langle\lambda_n\rangle_{n\leq N}$, a polytopological space $\langle X,\langle\mathcal T_\xi\rangle_{\xi<\Lambda}\rangle$ and a $\sf J$-frame $\langle W,\langle<_n\rangle_{n\leq N}\rangle$, we will say a map $f:X\to W$ is a {$\vec \lambda$-map} if it is a $\sf J$-map on $\langle X,\langle\mathcal T_{\lambda_n}\rangle_{n\leq N}\rangle$. These maps will allow us to focus on finitely many modalities at one time in the completeness proof.

Say a $\vec \lambda$-map is {\em suitable} if it is a surjective function of the form $f:(\Theta+1)\to W$, for the $0$-root $w_0$ of $W$ we have $f^{-1}(w_0)=\{\Theta\}$ and $\Theta$ lies in the range of $e$ (i.e., $\Theta=0$ or it is infinite and additively indecomposable). We use ${\rm hgt}(<_n)$ to denote the {\em height} of $<_n$, that is, the maximal $k$ such that there exist $w_0<_nw_1<_n\hdots<_nw_k$.

\begin{lemm}\label{main}
Given a $\sf J$-frame $\mathfrak W=\langle W,\langle{<_n}\rangle_{n\leq N}\rangle$ and ordinals $\vec\lambda=\langle\lambda_n\rangle_{n\leq N}$ all less than $\Lambda$, there exists an idyllic $\Lambda$-ambiance $\mathfrak X$ based on some $\Theta<e^{1+\Lambda}1$ and a suitable $\vec\lambda$-map $f:\mathfrak X\to \mathfrak W$.
\end{lemm}

\proof
We proceed as in the proof of an analogous result in \cite{BeklGabel:2011}.

Suppose $\langle W,\langle{<_n}\rangle_{n\leq N}\rangle$ is a $\sf J$-frame and $\lambda_0,\hdots\lambda_{N}$ are ordinals. We work by induction on $N$ with a secondary induction on ${\rm hgt}(<_0)$ to construct a suitable $\vec \lambda$-map. Withouth loss of generality, we assume $\lambda_0=0$, for otherwise we can always let $<_0=\varnothing$.

\paragraph{Case 1: $N=0$.} In this case it is known that there is an ordinal $\Theta=e\Theta'<\omega^\omega=e(\omega)$ and a suitable $d$-map $g:(\Theta+1)_{1}\to W$ (see \cite{glcomplete}). Note that since $\Theta<\omega^\omega$, it has no points of limit rank and hence the interval topology is already limit-maximal so that $\mathcal P(\Theta+1)$ is an idyllic algebra.

\paragraph{Case 2: ${\rm hgt}(<_0)=0$.} Here we have that $<_0=\varnothing$. For $0<n\leq N$ let $\lambda_n'=-\lambda_1+\lambda_n$ and consider the $\sf J$-frame $\langle W,\langle{<_{n+1}}\rangle_{0\leq n< N}\rangle$, where $<_0$ has been removed. By induction on $N$ we may assume there is an idyllic ambiance $\mathfrak X$ based on a polytopology $\vec{\mathcal T}$ on an ordinal $\Theta+1<e^{1+(-\lambda_1+\Lambda)}1$ and a suitable $\vec\lambda '$-map $g$ from $\mathfrak X$ onto $W$.

Let $\Omega=e^{1+\lambda_1} \le\Theta$ and use Lemma \ref{isbg} to construct a BG-polytopology
\[{\mathfrak Y}=\langle \Omega+1,\langle\mathcal S_\lambda\rangle_{\lambda<\lambda_1}\rangle.\]
Note that from $\Theta+1<e^{1+(-\lambda_1+\Lambda)}$ we obtain $\le\Theta<e^{-\lambda_1+\Lambda}$ and thus
\[e^{1+\lambda_1}\le\Theta<e^{1+\lambda_1} e^{-\lambda_1+\Lambda}1=e^{1+\lambda_1+ (-\lambda_1+\Lambda)}1=e^{1+\Lambda}1.\]

Then, by Theorem \ref{thereismaptheo} there is a $\Lambda$-reductive map \[f=\reductive^{\le\Theta}_{\lambda_1}:(\Omega+1)_{1+\lambda_1}\to(\Theta+1)_1,\]
so that by Lemma \ref{stilld} $f:{\mathfrak Y}^-_{\lambda_1}\to\mathfrak X_0$ is a $d$-map. By Lemma \ref{LemmBG}.\ref{LemmBGFour}, there exists a BG-space ${\mathfrak Y}_{\lambda_1}$ extending ${\mathfrak Y}^-_{\lambda_1}$ such that  $f:{\mathfrak Y}_{\lambda_1}\to {\mathfrak X}_0$ is a $d$-map.

Hence we may use Lemma \ref{polylift} to define BG-topologies $\langle\mathcal S_{\xi}\rangle_{\lambda_1<\xi\leq \lambda_N}$ on $\Omega+1 $ such that $f:{\mathfrak Y}_{\lambda_1+\zeta}\to{\mathfrak X}_{\zeta}$ is a $d$-map for all $\zeta$; in particular, for each $n\in (0,N]$ we have that $f:{\mathfrak Y}_{\lambda_n}\to{\mathfrak X}_{\lambda'_n}$ is a $d$-map.

It follows that $f:{\mathfrak Y}\to{\mathfrak X}$ is a $d$-lift, and by Lemma \ref{pullalg}, $f^{-1}\mathcal A$ is an idyllic $d$-algebra. Further, by Lemma \ref{dj} we know that $gf:{\mathfrak Y}\to W$ is a suitable $\vec\lambda$-map, as needed.

\paragraph{Case 3: ${\rm hgt}(<_0)=m>0$.} Let $w$ be the $0$-root of $W$ and $w_0,\hdots w_K$ be its $<_0$-daughters which are hereditary $1$-roots.

Let $V=\{w\}\cup \mathop\ll_1(w)$ and $W_k=\mathop\ll_0(w_k)$. Then we have, as in Case 2, an idyllic ambiance $\mathfrak X$ based on an ordinal $\Theta+1<e^{1+\Lambda}1$ and by induction on $N$ a $\vec{\lambda}$-map $f$ from $\Theta+1$ to $V$, as well as for each $k\leq K$ an idyllic ambiance $\mathfrak Y_k$ based on an ordinal $\Xi_k+1<e^{1+\Lambda}1$ and a suitable $\vec\lambda$-map $f_k$ from $\Xi_k+1$ onto $W_k$.

Let $\Xi=(\Xi_0+1)+(\Xi_1+1)+\hdots+\Xi_K$, $\mathfrak Y=\bigoplus_{k\leq K}\mathfrak Y_k$ and $\mathfrak Z=\mathfrak Y\otimes_d\mathfrak X$ with associated maps $\pi_0$ and $\pi_1$. Let $G_0$ be the domain of $\pi_0$ and $G_1$ be its complement. Define $g:\mathfrak Y\to W$ by \[g\big((\Xi_0+1)+\hdots+(\Xi_{k-1}+1)+\xi\big)=f_k(\xi)\]
and $h:\mathfrak Z\to W$ by
\[h(\xi)=
\begin{cases}
g \pi_0 (\xi)&\text{ if $\xi\in G_0$,}\\
f\pi_1 (\xi)&\text{ otherwise.}
\end{cases}
\]

Note that $\Lambda>0$ so that $e^{1+\Lambda}1$ is closed under sums and products and thus $\Xi\otimes_d\Theta<e^{1+\Lambda}1.$ Further, by Lemma \ref{prodil}, $\mathfrak Z$ is an idyllic ambiance.

Now, let us check that $h$ satisfies the conditions of Definition \ref{jmap}.
\paragraph{1.} We know that $\lambda_N>0$, since $N>0$. Thus $G_0,G_1$ are both $\lambda_N$-clopen and hence it is enough to check that $h:(\mathfrak Z_{\lambda_n}\upharpoonright G_j)\to \langle W,<_K\rangle$ is a $d$-map for $j=0,1$. But this is immediate from the assumption that $f,g$ are $d$-maps, as are the respective projections.

\paragraph{2.} Let $\lambda=\lambda_n$, $U$ be $\lambda$-open and $v\in h(U)$. We have that
\[h(U)=g\pi_0(G_0\cap U)\cup f\pi_1 (G_1\cap U).\]
First note that $G_0\cap U$ is $\lambda$-open since $G_0$ is $0$-open, hence $g\pi_0(G_0\cap U)$ is also $\lambda$-open given that $g\pi_0$ is a composition of $\lambda$-open maps. Meanwhile, for $\lambda>1$, $G_1$ is also $\lambda$-open from which it follows that $f\pi_1 (G_1\cap U)$ is $\lambda$-open as well.

It remains to check that $h(U)$ contains a $0$-neighborhood around any $v\in f\pi_1 (G_1\cap U)$. So suppose $u<_0v$. Since $g$ is onto $\bigcup W_n$, there is some $\delta\leq \Xi$ such that $g(\delta)=u$. By Lemma \ref{LemmDprodTop}.\ref{LemmDprodTop2} and using the fact that $G_1\cap U\not=\varnothing$ (otherwise $v$ would not exist), there is $\gamma\in U\cap \pi_0^{-1}\delta$. It follows that $h(\gamma)=g\pi_0(\gamma)=u$, as desired.

\paragraph{3.}

If $n>0$ and $v$ is an $n+1$ root the claim follows form the assumption that $f,g$ were already $\sf J$-maps and the respective projections are $d$-maps. Meanwhile, if $n=0$ and $v<_0w$, we may use the fact that $g$ is a $\langle\lambda_n\rangle_{1\leq n\leq N}$-map, since here it follows that $h^{-1}(\mathop\ll_0(v))\subseteq G_0$ so that
\[h^{-1}(\mathop{\ll_0}(v))=(g\pi_0)^{-1}(\mathop{\ll_0}(v)).\]
But $g^{-1}(\mathop \ll_0(v))$ is a $0$-open subset of $\Xi$, hence $(g\pi_0)^{-1}(\mathop \ll_0(v))$ is open in $\mathfrak Z$. Similarly for $\{v\}\cup\mathop\ll_0(v)$.

If $v\ll_1 w$, then $\mathop{<_0}(v)=\mathop{<_0}(w)$, and thus $h^{-1}(\mathop{<_0}(v))=G_0$. But then we see that
\[h^{-1}(\mathop{\ll_0}(v))=G_0\cup \pi_1^{-1}f^{-1}(\mathop{\ll_0}(v));\]
since by assumption $f$ is a $\sf J$-map, $f^{-1}(\mathop{\ll_0}(v))$ is $0$-open in $\mathfrak Z\upharpoonright G_1$, and therefore by Lemma \ref{LemmDprodTop}.\ref{LemmDprodTop3}, $G_0\cup \pi_1^{-1}f^{-1}(\mathop{\ll_0}(v))$ is $0$-open in the $d$-product $\mathfrak Z$. The argument for $\{v\}\cup\mathop\ll_0(v)$ is analogous.

Finally, note that $f^{-1}(w)=\{\Xi\otimes_d\Theta\}$ so that $f^{-1}(\mathop \ll_0(w))=W\setminus \{w\}=[0,\Xi\otimes_d\Theta)$, which is open, as is $f^{-1}(\{w\}\cup\mathop \ll_0(w))=[0,\Xi\otimes_d\Theta]$.

\paragraph{4.}

For $w$ we have that $h^{-1}(w)=\{\Xi\otimes_d\Theta\}$, which is obviously $\lambda_n$-discrete for any $n$. If $v<_0w$ then $h^{-1}(v)=(g\pi_0)^{-1}(v)$ which is discrete, as $g\pi_0$ is a $\sf J$-map. 

If $v\ll_1w$ then $h^{-1}(v)=(f\pi_1)^{-1}(v)$, which similarly must be discrete.
\endproof

With this, we are ready to state and prove our main result.

\begin{theorem}\label{bigone}
${\sf GLP}_\Lambda$ is complete for the class of idyllic ambiances based on some $\Theta< e^{1+\Lambda}1$.
\end{theorem}

\proof
Suppose that $\phi$ is consistent over $\glp_\Lambda$. Then, by Lemma \ref{condenselemm}, $\phi^{\sf c}$ is consistent over $\glp_N$. By Lemma \ref{glptoj}, $M^+(\phi)\wedge\phi$ is consistent over $\sf J$, and thus by Lemma \ref{jcomp}, we have a tree-like  $\sf J$-model $\langle W,\langle <_n\rangle_{n\leq N},\lb\cdot\rb\rangle$ satisfying $M^+(\phi)\wedge\phi$. We may then use Lemma \ref{main} to find an idyllic ambiance $\mathfrak X$ based on an ordinal $\Theta<e^\Lambda\omega$ and a surjective $\vec\lambda$-map $f:{\mathfrak X}\to  W$.

Then, by Lemma \ref{truthpres}, there is a valuation on $\Theta$ agreeing with $f^{-1}\lb\cdot\rb$ on ${\rm sub}(\phi)$. Since $f$ is surjective, $f^{-1}\lb\phi\rb\not=\varnothing$, and thus $\phi$ is satisfed on $\mathfrak X$, as desired.
\endproof

As corollaries we get a sequence of completeness results:

\begin{cor}
${\sf GLP}_\Lambda$ is complete for both the class of BG-spaces and the class of shifted Icard ambiances based on $e^{1+\Lambda}1$.

Further, the variable-free fragment ${\sf GLP}^0_\Lambda$ is complete for the class of simple Icard ambiances based on $e^{\Lambda}1$.
\end{cor}

\proof
Completeness for BG-spaces and shifted Icard ambiances is immediate from Theorem \ref{bigone}, as idyllic ambiances may be seen as either kind of structure.

Meanwhile, given any Icard ambiance based on an algebra $\mathcal A$ and satisfying a closed formula $\phi$, we use Lemma \ref{issimp} to note that all valuations of $\phi$ and its subformulas are simple, and hence we obtain a simple ambiance satisfying $\phi$ by replacing $\mathcal A$ by the class of simple sets.

This gives us an ambiance on $\widehat{{\mathfrak {Ic}}}^{e^{1+\Lambda}1}_\Lambda$, the shifted Icard space. To pass to an ambiance based on ${\mathfrak {Ic}}^{e^{\Lambda}1}_\Lambda$, note that $\le:\widehat{{\mathfrak {Ic}}}^{e^{1+\Lambda}1}_\Lambda\to{\mathfrak {Ic}}^{e^{\Lambda}1}_\Lambda$ is a $d$-map and thus by an easy induction preserves valuations of closed formulas.
\endproof

The completeness result for simple ambiances is not new, as it was already proven by Icard for ${\sf GLP}_\omega$ in \cite{Icard:2009} and by Joosten and I for arbitrary ${\sf GLP}_\Lambda$ in \cite{glpmodels}. However, the current argument is quite different from those used in previous works.

\section{Worms and the lower bound}\label{lowsec}

As it turns out, our bound of $e^{1+\Lambda}1$ is sharp. To show this, let us consider {\em worms}.

A worm is a formula of the form \[\langle\lambda_0\rangle\hdots\langle\lambda_I\rangle\top.\]
These formulas correspond to iterated consistency statements, and indeed can be used to study the proof-theoretic strength of many theories related to Peano Arithmetic, as Beklemishev has shown \cite{Beklemishev:2004}.

Worms are well-ordered by their {\em consistency strength}. Let us denote the set of worms with entries less than $\Lambda$ by ${\mathbb W}^\Lambda$; then, given worms $\mathfrak v,\mathfrak w\in {\mathbb W}^\Lambda$, define $\mathfrak v\lhd\mathfrak w$ if $\glp_\Lambda\vdash\mathfrak w\to\ps\mathfrak v$.

The relation $\lhd$ we have just defined is a well-order \cite{Beklemishev:2004,WellOrdersI}. Thus we may compute the order-type of a worm $\mathfrak w\in\mathbb W^\Lambda$:

\[o(\mathfrak w)=\sup_{\mathfrak v\lhd\mathfrak w}(o(\mathfrak v)+1).\]

\begin{lemm}
Let $\mathfrak w$ be a worm, $\mathfrak X=\langle X,\langle\mathcal T_\lambda\rangle_{\lambda<\Lambda},\lb\cdot\rb\rangle$ be a GLP-model and $x\in X$.

Then, if $x\in \lb\mathfrak w\rb$, it follows that $\uprho(x)\geq o(\mathfrak w)$.
\end{lemm}

\proof
By induction on $o(\mathfrak w)$. For the base case, note that if $o(\mathfrak w)=0$ then we vacuously have that $x\in\lb\mathfrak w\rb$ implies $\uprho(x)\geq 0$.

For the inductive step, if $\mathfrak v\lhd\mathfrak w$ and $U$ is any neighborhood of $x$, since $\vdash \mathfrak w\to\ps\mathfrak v$, it follows that there is $y\in U$ satisfying $\mathfrak v$. By induction on $\mathfrak v\lhd\mathfrak w$, we have that $\uprho(y)\geq o(\mathfrak v)$.

We then see that \[\uprho(x)\geq\displaystyle \sup_{\mathfrak v\lhd\mathfrak w}(o(\mathfrak v)+1)=o(\mathfrak w).\]
\endproof

It will be convenient to review the calculus for computing $o$ that is given in \cite{WellOrdersI}. First, if $\mathfrak v=\langle\xi_1\rangle\hdots\langle\xi_N\rangle\top$ and $\mathfrak w=\langle\zeta_1\rangle\hdots\langle\zeta_M\rangle\top$, define
\[\mathfrak v\ps\mathfrak w=\langle\xi_1\rangle\hdots\langle\xi_N\rangle\langle 0\rangle\langle\zeta_1\rangle\hdots\langle\zeta_M\rangle\top.\]

Further, if $\alpha$ is any ordinal, set \[\alpha\uparrow\mathfrak w=\langle\alpha+\zeta_1\rangle\hdots\langle\alpha+\zeta_M\rangle\top.\]

\begin{lemm}\label{wormlemm}
Let $\mathfrak v,\mathfrak w$ be worms and $\alpha$ an ordinal.

Then,
\begin{align}
o(\top)&=0\\
o(\mathfrak v\ps\mathfrak w)&=o(\mathfrak w)+1+o(\mathfrak v)\label{second}\\
o(\alpha\uparrow\mathfrak w)&=e^\alpha o(\mathfrak w).
\end{align}

\end{lemm}

\begin{examp}
Consider $\mathfrak w=\langle \omega+1\rangle \langle \omega\rangle\langle \omega+1\rangle\top$; let us compute $o(\mathfrak w)$. Note that $\mathfrak w=\omega\uparrow\mathfrak v$, where $\mathfrak v=\langle 1\rangle \langle 0\rangle\langle 1\rangle\top$, so that $o(\mathfrak w)=e^\omega o(\mathfrak v)$. Now, $o(\mathfrak v)=o(\langle 1\rangle\top)+1+o(\langle 1\rangle\top)$, and we have that $o(\langle 1\rangle\top)=o(1\uparrow \langle 0\rangle\top)= eo(\top\ps\top)=e(0+1+0)=\omega$, so that $o(\mathfrak v)=\omega+1+\omega=\omega+\omega$ and $o(\mathfrak w)=e^\omega(\omega+\omega)=\varepsilon_{\omega+\omega}$.
\end{examp}

\begin{theorem}
${\sf GLP}_\Lambda$ is incomplete for the class of BG-spaces or Icard ambiances based on any fixed $\Theta<e^{1+\Lambda} 1$.
\end{theorem}

\proof
Note that $\le(e^{1+\Lambda} 1)=e^\Lambda 1$. Now, if $\glp_\Lambda$ is complete for the class of models based on $\Theta$, in particular any worm $\mathfrak w\in {\mathbb W}^\Lambda$ must be satisfiable on one such model $\mathfrak X$, which implies by Lemma \ref{wormlemm} that there must be $\vartheta\in \Theta$ with $ \le(\vartheta)\geq o(\mathfrak w)$. Thus it suffices to show that
\[\sup_{\mathfrak w\in{\mathbb W}^\Lambda}o(\mathfrak w)\geq e^\Lambda 1.\]

To do this, first assume $\Lambda=\lambda+1$. Then we have that
\[o(\langle\lambda\rangle^n\top) =e^\lambda o(\langle 0\rangle^n\top)=e^\lambda n.\]

The last equality is obtained by repeated applications of (\ref{second}).

But then, by Proposition \ref{prophyp}.\ref{prophypthree} we see that
\[\uprho(\mathfrak X)\geq\lim_{n\to\omega}e^\lambda n=e^\lambda \omega=e^\lambda e 1=e^{\lambda+1}1.\]

Meanwhile, if $\Lambda \in\sf Lim$, \[\uprho(\mathfrak X)\geq\displaystyle\sup_{\lambda<\Lambda}o(\langle\lambda\rangle\top) =\displaystyle\sup_{\lambda<\Lambda}e^\lambda 1=e^\Lambda 1.\]

In either case $\uprho({\mathfrak X})=\sup_{\vartheta<\Theta}\le\vartheta\geq e^\Lambda 1$, from which it follows by Lemma \ref{leftinv} that $\Theta\geq e^{1+\Lambda}1$.
\endproof

\section{Concluding remarks}

The goal of this paper was essentially to answer two main questions. The first is perhaps not so much {\em Is $\glp_\Lambda$ complete for its topological semantics independently of $\Lambda$?} as, rather, {\em What is needed to construct topological models of $\glp_\Lambda$?} For this we had to introduce several tools that were not required in the case $\Lambda=\omega$. Most notable is the use of hyperlogarithms and -exponents, already employed in \cite{glpmodels} to study models of the closed fragment, and the addition of new $d$-maps to our toolkit. Aside from possible connections to proof theory, these are novel constructions in scattered topology and might spark some independent interest.

The second question is, {\em Are there good constructive semantics for $\glp_\Lambda$?} Icard ambiances are a possible answer to this question. Not only are the topologies easily definable, unlike the non-constructive BG-topologies, but if one analyzes the proof of Lemma \ref{main}, all sets that appear in valuations are constructive as well. As such, Icard ambiances may be well-suited for applications in the proof theory of systems much stronger than Peano Arithmetic -- perhaps the ultimate motivation for contemporary work in provability logic.

\bibliographystyle{plain}

\end{document}